\documentclass[a4paper,12pt,reqno]{amsart}

\usepackage{amssymb}
\usepackage[english]{babel}
\usepackage{a4wide}
\usepackage{amsmath,amssymb}
\usepackage{algorithm}
\usepackage{graphicx}
\usepackage{subfig}
\usepackage{hyperref}
\usepackage{psfrag,float}
\usepackage{graphicx}

\usepackage{amsthm}
\newtheorem{thm}{Theorem}[section]
\newtheorem{rem}[thm]{Remark}

\DeclareMathOperator{\cof}{cof}
\DeclareMathOperator{\Div}{div}

\title{Numerical simulation of nonlinear continuity equations by evolving diffeomorphisms}

\author[label1]{Jos\'e A. Carrillo$^*$}
\author[label2]{Helene Ranetbauer$^\ddagger$}
\author[label3]{Marie-Therese Wolfram$^\S$}

\begin{document}
\maketitle
\renewcommand{\thefootnote}{\fnsymbol{footnote}}
\footnotetext[1]{Department of Mathematics, Imperial College London, London SW7 2AZ, UK;\\ email: carrillo@imperial.ac.uk}
\footnotetext[3]{Radon Institute for Computational and Applied Mathematics, Austrian Academy of Sciences, Altenberger Strasse 69, 4040 Linz, Austria;\\ email: helene.ranetbauer@ricam.oeaw.ac.at}
\footnotetext[4]{Mathematics Institute, University of Warwick, Coventry CV4 7AL, UK and Radon Institute for Computational and Applied Mathematics, Austrian Academy of Sciences, Altenberger Strasse 69, 4040 Linz, Austria;\\ email: m.wolfram@warwick.ac.uk}
\renewcommand{\thefootnote}{\arabic{footnote}}
\begin{abstract}
  In this paper we present a numerical scheme for nonlinear continuity equations, which is based on the gradient flow formulation of an energy functional with respect to the quadratic transportation distance. It can be applied to a large class of nonlinear continuity equations, whose dynamics are driven by internal energies, given external potentials and/or interaction energies. The solver is based on its variational formulation as a gradient flow with respect to the Wasserstein distance. Positivity of solutions as well as energy decrease of the semi-discrete scheme are guaranteed by its construction. We illustrate this properties with various examples in spatial dimension one and two.
\end{abstract}


\markboth{Carrillo, Ranetbauer, Wolfram }{FEM for evolving diffeomorphisms}

\section{Introduction}
In this work we propose a numerical method for solving nonlinear continuity equations of the form:
\begin{subequations}\label{e:cont}
\begin{align}
\partial_t \rho &= -\nabla \cdot [ \rho v] := \nabla \cdot \left[ \rho \nabla \left(U'(\rho) + V + W *\rho\right)\right]\\
\rho(0,\cdot) &= \rho_0,
\end{align}
\end{subequations}
where  $\rho = \rho(t,x)$ denotes the unknown time dependent probability density and $\rho_0 = \rho_0(x)$ a given initial probability density on $\Omega \subseteq \mathbb{R}^d$. The function $U: \mathbb{R}^+ \rightarrow \mathbb{R}$ is an internal energy, $V=V(x): \mathbb{R}^d\rightarrow \mathbb{R}$ a given potential and $W = W(x): \mathbb{R}^d \rightarrow \mathbb{R}$ an interaction potential. Equation \eqref{e:cont} can be interpreted as a gradient flow with respect to the Euclidean Wasserstein distance of the free energy or entropy
\begin{align*}
\mathcal{E}(\rho) = \int_{\mathbb{R}^d} U(\rho(x)) dx + \int_{\mathbb{R}^d} V(x) \rho(x) dx + \frac{1}{2}\int_{\mathbb{R}^d\times \mathbb{R}^d} W(x-y) \rho(x) \rho(y) dx dy.
\end{align*}
Then the velocity field $v = v(t,x)$ corresponds to $v = -\nabla \frac{\delta \mathcal{E}}{\delta \rho}$, and the free energy is dissipated
along the trajectories of equation \eqref{e:cont}, i.e.
\begin{align*}%\label{e:energydissipation}
\frac{d}{dt} \mathcal{E}(\rho)(t) = -D(\rho) \equiv -\int_{\mathbb{R}^d} \lvert v(t,x) \rvert^2 \rho(t,x) dx.
\end{align*}
Here $D(\rho)$ denotes the so called entropy dissipation functional. 

The gradient flow formulation detailed above provides a natural framework to describe the evolution of densities and it has been successfully used to model transportation processes in the life and social sciences. Examples of it include the heat equation with $U(s) = s \ln s$, $V = W = 0$ or the porous medium and fast diffusion equation with $V = W = 0$ and $U(s) = s^{m} /(m-1)$ for $m>1$ and $0 < m < 1$ respectively. If the dynamics are in addition driven by a given potential $V$, we refer to \eqref{e:cont} as a Fokker-Planck type equation.
Gradient flow techniques have also been applied to aggregation equations, which correspond to \eqref{e:cont} for a given interaction potential $W$ and $U=V=0$.
These models appear naturally to describe spatial shapes of collective dynamics in general and in particular for animal swarms, for example the motion of bird flocks or fish schools, cf. \cite{mogilner1999non,TP,TPL,d2006self}. See also the reviews \cite{review,review2} and the references therein. Applications in physics include the field of granular media \cite{carrillo2003kinetic,T2000} or material sciences \cite{HP2005}. A common assumption in collective dynamics is that $W$ is a radial function, i.e. $W = W(\lvert x\rvert)$. Hence interactions among individuals depend on their distance only. The interaction dynamics are often driven by attractive and repulsive forces, mimicking the tendency of individuals to stay close to
the group but maintain a minimal distance. A popular choice among attractive-repulsive potentials is the Morse potential, that is
\begin{align*}%\label{e:morse}
W(r) = -C_A e^{-r/l_A} + C_R e^{-r/l_R},
\end{align*}
where $C_A$ and $C_R$ are the attractive and repulsive strength and $l_A,~l_R$ their respective length scales, see \cite{LTB,BT,CMP,CHM,Bertozzietal,ABCV}. Also power-law potentials of the form
\begin{align*}%\label{e:powerlaw}
W(r) = \frac{\lvert r \rvert^a}{a} - \frac{\lvert r \rvert^b}{b}, ~a>b,
\end{align*}
have been used thoroughly. We refer to the works of \cite{Bertozzietal,BCLR,BCLR2,fellner2010stable,fellner2011stability} for more details. Purely attractive potentials are often of the form $W(r) = \lvert r \rvert^{a},~a>0$. In this case the density of particles collapses in finite or in infinite time and converges to a Delta Dirac located at the center of mass for certain range of values, see \cite{BertozziCarilloLaurent,BH1,BH2,CDFLS2011}. 
Finally, aggregation-diffusion equations, in which both the linear or nonlinear diffusion term modelling repulsion and the aggregation term modeling attraction are present, are ubiquitous models in physics and biology, including all different versions of the parabolic-elliptic Keller-Segel model for chemotaxis, see for instance \cite{HV,V,blanchet2008infinite,blanchet2008convergence,BCL,YaoBer} and the references therein.

The gradient flow structure of \eqref{e:cont} provides a natural framework to construct solutions using variational schemes, usually called the JKO algorithm, as proposed in the seminal works of Jordan, Kinderlehrer and Otto \cite{jordan1998variational,O2001}. The variational formulation provides an underlying structure for the construction of numerical methods, which have inherent advantages such as build-in positivity and free-energy decrease.  There has been an increasing interest in the development of such methods in the last years, for example variational Lagrangian schemes such as \cite{gosse2006lagrangian,gosse2006identification,WW2010,during2010gradient, blanchet2008convergence,MO2014,carrillo2015numerical, junge2015fully} or finite volume schemes as in \cite{CCH2015}. The common challenge of all methods is the high computational complexity, often restricting them to spatial dimension one. They involve for example the solution of the Monge-Kantorovich transportation problem between measures, which require the computation of the Wasserstein distance. In 1D the Wasserstein distance corresponds to the $L^2$ distance between inverses of the distribution functions. In higher space dimensions its computation involves the solution of an optimal control problem itself, which poses a significant challenge for the development of numerical schemes. There are few results on the numerical analysis of these schemes, for example Matthes and co-workers provided first results on the convergence in 1D in \cite{MO2014}. Recent works \cite{benamou2014discretization} are making use of recent developments in the fast computation of optimal transportation maps to discretize the JKO steps in the variational scheme. 

Many of the Lagrangian methods are based on the ``discretize-then-minimize'' strategy, see for example \cite{WW2010,during2010gradient,MO2014,benamou2014discretization,carrillo2015numerical,junge2015fully}. However, we follow closely the advantage of the 1D formulation in terms of the pseudo-inverse function, see also \cite{gosse2006lagrangian,gosse2006identification, blanchet2008convergence}, which corresponds to a ``minimize-then-discretize'' strategy. In higher dimensions, this strategy is based on the variational formulation of \eqref{e:cont} for a diffeomorphism mapping the uniform density to the unknown density $\rho$. Evans et al. proved in \cite{ESG2005} that this approach corresponds to a class of $L^2$-gradient flows of functionals on diffeomorphisms. This equivalent formulation served as a basis for the numerical solver of Carrillo and Moll \cite{carrillo2009numerical}, which uses the variational formulation and an explicit in time as well as a finite difference discretization in space. We shall follow their approach but propose an implicit time stepping and a spatial discretization, which is based on finite differences in 1D and finite elements in 2D. We would like to mention that the our scheme coincides with the implicit in time finite difference discretization of the pseudo-inverse function proposed by Blanchet et al. for the 1D Patlak-Keller-Segel model, see \cite{blanchet2008convergence}. The finite element discretization allows us to consider more general computational domains as well as triangular meshes. This is advantageous in the case of radially symmetric solutions, which develop asymptotically for large times for many types of aggregation-diffusion equations. Due to the implicit in time stepping no CFL condition as in \cite{carrillo2009numerical} is necessary. Moreover, if our dicretization were convergent for the implicit semidiscretization of the PDE system satisfied by the diffeomorphisms, the numerical scheme would be convergent for an exact JKO step of the the variational scheme.

This paper is organized as follows: we review the underlying ideas of the proposed variational scheme in Section \ref{s:gradflow}. In Section \ref{s:diffeo}
we discuss the numerical algorithm; Section \ref{s:preandpost} details the necessary pre- and postprocessing steps. Finally we present extensive numerical simulations for a large class of nonlinear aggregation-diffusion equations in Section \ref{s:numerics} including blow-up profiles and complicated asymptotic behaviors showing the flexibility of the method to cope both with diffusive and aggregation behavior in the same model. Section \ref{s:conclusion} is devoted to discuss the main aspects of our approach together with future perspectives.

%%%%%%%%%%%%%%%%%%%%%%%%%%%%%%%%

\section{Gradient flow formulation}\label{s:gradflow}

We start by briefly reviewing the variational formulation of \eqref{e:cont} as a gradient flow with respect to a particular distance and energy.  Let us consider the quadratic
Wasserstein distance between two probability measures $\mu \in \mathcal{P}(\mathbb{R}^d)$ and $\nu \in \mathcal{P}(\mathbb{R}^d)$ given by
\begin{align*}%\label{e:wasserstein}
d_W^2(\mu, \nu) := \inf_{T: \nu = T \# \mu} \int_{\mathbb{R}^d} \lvert x-T(x) \rvert^2 d\mu(x)\,,
\end{align*}
for any $\mu$ absolutely continuous with respect to Lebesgue, see \cite{V2003} for a general definition and its properties.
It is well known that solutions of \eqref{e:cont} can be constructed via the so-called JKO scheme, see \cite{jordan1998variational}, which corresponds to solving
\begin{align}\label{e:varscheme}
\rho^{n+1}_{\Delta t} \in \arg \min_{\rho \in \mathcal{K}}\left\{ \frac{1}{2 \Delta t} d_W^2(\rho^n_{\Delta t}, \rho) + \mathcal{E}(\rho)\right\},
\end{align}
for a fixed time step $\Delta t > 0$ and $\mathcal{K} = \lbrace \rho \in L^1_+(\mathbb{R}^d) : \int_{\mathbb{R}^d} \rho(x) dx = M, \lvert x \rvert^2 \rho \in L^1(\mathbb{R}^d) \rbrace$.
Hence \eqref{e:varscheme} can be understood as a time discretization of an abstract gradient flow equation in the space of probability measures. It has been proven that
solutions of \eqref{e:varscheme} converge to the solutions of \eqref{e:cont} first order in time, see \cite{jordan1998variational, O2001, V2003, AGS2005} for more detailed results on the analysis.

Evans et al. \cite{ESG2005} showed  that there is a connection between the theory of steepest descent schemes with respect to the Euclidean transport distance and the $L^2$-gradient
flows of polyconvex functionals of diffeomorphisms. Since this formulation serves as a basis of the proposed numerical solver, we will review the main results in the following.
Let $\tilde{\Omega}$ be a smooth, open, bounded and connected subset of $\mathbb{R}^d$ and $\Omega$ be an open subset of $\mathbb{R}^d$. 
Let $\mathcal{D}$ denote the set of diffeomorphisms from $\bar{\tilde{\Omega}}$ to $\bar{\Omega}$, mapping $\partial \tilde{\Omega}$ onto $\partial \Omega$. 
Furthermore we consider the energy functional:
\begin{align*}
\mathcal{I}(\Phi) = \int_{\tilde{\Omega}} \Psi(\det D\Phi) dx + \int_{\tilde{\Omega}} V(\Phi(x)) dx + \frac{1}{2} \int_{\tilde{\Omega} \times \tilde{\Omega}} W(\Phi(x)-\Phi(y)) dx dy,
\end{align*}
where $\Psi(s)=sU\left(\frac{1}{s}\right)$.
Ambrosio et al \cite{ambrosio2006stability} clarified even further this subtle relation found in \cite{ESG2005}, by showing that the $L^2$-gradient flow of $\mathcal{I}(\Phi)$ with $V=W=0$ given by
\begin{align}\label{e:varl2}
\Phi^{n+1}_{\Delta t} \in \arg \min_{\Phi \in \mathcal{D}} \left\{ \frac{1}{2 \Delta t} \lVert \Phi^n_{\Delta t} - \Phi \rVert_{L^2(\Omega)} + \mathcal{I}(\Phi) \right\}
\end{align} 
is well defined and converges to the solutions of the nonlinear PDE system
\begin{align*}
\frac{\partial \Phi}{\partial t} &=  \nabla \cdot \left[ \Psi'(\det D\Phi) (\cof D\Phi)^T \right] ,
\end{align*}
where $D\Phi$ is the Jacobian matrix of $\Phi$ and $\cof D\Phi$ the corresponding cofactor matrix. This connection was generalized subsequently in \cite{carrillo2009numerical} to the case of interaction and potential energies, leading to the nonlinear PDE system
\begin{align}\label{e:diffeo}
\frac{\partial \Phi}{\partial t} &=  \nabla \cdot \left[ \Psi'(\det D\Phi) (\cof D\Phi)^T \right] - \nabla V \circ \Phi - \int_{\tilde{\Omega}} \nabla W(\Phi(x)- \Phi(y)) dy,
\end{align}\\
The diffeomorphism $\Phi$ maps a given reference density, for instance the uniform density on a reference domain $\tilde{\Omega}$, to the unknown density $\rho$ in $\Omega$. Therefore the density $\rho \in \mathcal{K}$ can be calculated via
$\rho = \Phi \# \mathcal{L}^N$ for every diffeomorphism $\Phi \in \mathcal{D}$, assuming that $|\tilde{\Omega}|=M$, where $\mathcal{L}^N$ denotes the $N$-dimensional Lebesgue measure on the reference domain $\tilde{\Omega}$. An equivalent formulation in the case of a sufficiently smooth diffeomorphism $\Phi$ is given by
\begin{align}\label{e:varchange}
\rho(\Phi(x)) \det(D \Phi(x)) = 1.
\end{align}
Hence we can interpret equation \eqref{e:diffeo} as the Lagrangian representation of the original nonlinear continuity equation \eqref{e:cont} in Eulerian coordinates.

%%%%%%%%%%%%%%%%%%%%%%%%%%%%%%%%%%%%%%%%

\section{Spatial and temporal discretization of the Lagrangian representation}\label{s:diffeo}

\noindent In this section we present the details of the spatial and temporal discretization of equation \eqref{e:diffeo}, focusing on the implicit in time scheme and the
subsequent finite dimensional approximation of the nonlinear semi-discrete equations. The latter is discretized by a finite difference scheme in 1D and a finite element 
method in 2D. In the following we shall only present the finite element discretization of \eqref{e:diffeo}, since the 1D finite difference scheme is detailed
in \cite{blanchet2008convergence,carrillo2009numerical} already. 

Let $\Delta t$ denote the discrete time step, $t^{n+1} = (n+1) \Delta t$ and $\Phi^{n+1}$ the solution $\Phi = \Phi(t,x)$ at time $t^{n+1}$. Then the implicit in time discretization of \eqref{e:diffeo} reads as
\begin{align}\label{e:implicitdiffeo}
\frac{\Phi^{n+1}-\Phi^n}{\Delta t} = \,& \nabla \cdot [\Psi'(\det D\Phi^{n+1}) (\cof D\Phi^{n+1})]\nonumber \\&-   \nabla V(\Phi^{n+1}) - \int_{\tilde{\Omega} } \nabla W(\Phi^{n+1}(x)-\Phi^{n+1}(y)) dy.
\end{align}
Let us consider test functions $\varphi = \varphi(x) \in H^1(\tilde{\Omega})$. Then the nonlinear operator $F$ defined via \eqref{e:implicitdiffeo} is given by:
\begin{align}\label{e:F}
\begin{split}
F(\Phi,\varphi) = \frac{1}{\Delta t} &\int_{\tilde{\Omega}}(\Phi^{n+1}-\Phi^n) \varphi(x) dx + \int_{\tilde{\Omega}} \Psi'(\det D\Phi^{n+1}) (\cof D\Phi^{n+1}) \nabla \varphi(x) dx \\
+ &\int_{\tilde{\Omega}} \nabla V(\Phi^{n+1}) \varphi(x) dx +\int_{\tilde{\Omega}} \left[\int_{\tilde{\Omega}} \nabla W(\Phi^{n+1}(x)-\Phi^{n+1}(y)) dy \right]  \varphi(x) dx .
\end{split}
\end{align}
We use lowest order $H^1$ conforming finite elements, also known as hat functions, for the discretization of $\Phi$, more precisely piece-wise linear functions. The discrete diffeomorphism $\Phi_h = \Phi_h(x_1, x_2)$ can be written as 
\begin{align*}
\Phi_h(x_1,x_2) = \sum_j \begin{pmatrix}\Phi^1_j\\[2mm] \Phi^2_j\end{pmatrix} \varphi_j(x_1,x_2),
\end{align*}
where the index $j$ corresponds to the nodal degrees of freedom.
We solve the nonlinear operator equation $F(\Phi,\varphi)=0$ using a Newton Raphson method in every time step (and drop the subscript $h$ to enhance readability in the following). To do so, we calculate the Jacobian matrix $DF$ of \eqref{e:F} and determine the Newton update $\Upsilon^{n+1,k+1}$ via
\begin{align*}
DF(\Phi^{n+1,k}, \varphi) \Upsilon^{n+1,k+1} =-F(\Phi^{n+1,k},\varphi),
\end{align*}
for all test functions $\varphi(x) \in H^1(\tilde{\Omega})$. The index $k$ corresponds to the Newton iteration and $n$ to the temporal discretization. Note that the Jacobian matrix $\Upsilon^{n+1,k+1}$ is a full matrix and has no sparse structure due to the convolution operator $W$. The Newton updates are calculated via
\begin{align*}
\Phi^{n+1,k+1} = \Phi^{n+1,k} + \alpha \Upsilon^{n+1,k+1},
\end{align*} 
where $\alpha$ is a suitable damping parameter. The Newton iteration is terminated if
\begin{align*}
  \lvert F(\Phi^{n+1,k+1},\varphi) \rvert \leq \epsilon_1 \text{ or } \lVert \Phi^{n+1,k+1}-\Phi^{n+1,k} \rVert \leq \epsilon_2,
\end{align*} 
with given error bounds $\epsilon_1$ and $\epsilon_2$. Note that due to the implicit in time discretization no CFL type condition for the time step $\Delta t$, as in \cite{carrillo2009numerical}, is necessary.

The presented $L^2$-gradient flow \eqref{e:varl2} is only valid if the total mass is conserved. This can be implemented by adopting the image domain in 1D or by considering equation \eqref{e:cont} with no flux boundary conditions in 2D, i.e.
\begin{align*}
v \cdot n = 0 \text{ on } \partial \Omega,
\end{align*}
where $n$ denotes the outer unit normal vector on the boundary $\partial \Omega$ and $v$ is given by \eqref{e:cont}. 
We shall only consider diffeomorphisms which map $\partial \tilde{\Omega}$ onto $\partial \Omega$ without rotations as in \cite{carrillo2009numerical}.
Then the corresponding natural boundary conditions for the equation in Lagrangian formulation are given by
\begin{align}\label{e:bc}
  n^T (\cof D \Phi)^T \frac{\partial \Phi}{\partial t} = (\cof D\Phi) n \cdot \frac{\partial \Phi}{\partial t} = 0.
\end{align} 
\noindent Note that the boundary conditions \eqref{e:bc} have to be checked separately for every computational domain, see for example \cite{carrillo2009numerical} for the
discussion of appropriate boundary conditions in the case of the unit square.  In the case of a circle with radius $R$ we have:
\begin{align*}
\frac{1}{R} \begin{pmatrix} x_1 \\ x_2 \end{pmatrix}
\begin{pmatrix} \frac{\partial \Phi_2}{\partial x_2} &-\frac{\partial \Phi_1}{\partial x_2} \\ -\frac{\partial \Phi_2}{\partial x_1} &\frac{\partial \Phi_1}{\partial x_1} \end{pmatrix} 
\begin{pmatrix} \frac{\partial \Phi_1}{\partial t} \\ \frac{\partial \Phi_2}{\partial t} \end{pmatrix} &= 0. 
\end{align*}
Using radial coordinates $x_1 = R \cos \theta$ and $x_2 = R \sin \theta$, we obtain:
\begin{align*}
\cos \theta (\frac{\partial \Phi_1}{\partial t} \frac{\partial \Phi_2}{\partial x_2} - \frac{\partial \Phi_2}{\partial t} \frac{\partial \Phi_1}{\partial x_2}) + \sin \theta(-\frac{\partial \Phi_1}{\partial t}
\frac{\partial \Phi_2}{\partial x_1} + \frac{\partial \Phi_2}{\partial t}\frac{\partial \Phi_1}{\partial x_1}) &= 0.
\end{align*} 
Since $\frac{\partial \Phi_1}{\partial \theta} = -\sin \theta \frac{\partial \Phi_1}{\partial x} = \cos \theta \frac{\partial \Phi_1}{\partial x_2}$, we conclude that
\begin{align}\label{e:bc2}
\sin \theta \frac{\partial \Phi_2}{\partial t} \frac{\partial \Phi_1}{\partial x_1}+ \cos \theta \frac{\partial \Phi_1}{\partial t}\frac{\partial \Phi_2}{\partial x_2} &= 0.
\end{align}

Hence, in the case of a circle equation \eqref{e:bc2} is equivalent to \eqref{e:bc}.\\
If we assume that the boundary of $\tilde{\Omega}$ is mapped onto the boundary of $\Omega$, the corresponding diffeomorphism is given by
\begin{align*}
\Phi_1(t,x) = \Phi_2(t,x) = Id,
\end{align*}
which implies that $\frac{\partial \Phi_1}{\partial t} = \frac{\partial \Phi_2}{\partial t} = 0$ on the boundary. This choice ensures that the equivalent boundary conditions given by \eqref{e:bc2} are satisfied. \\

%%%%%%%%%%%%%%%%%%%%%%%%%%%%%%%%%%%%%%%%%%

\section{Pre-processing and post-processing: calculating the initial diffeomorphism  and the final density}\label{s:preandpost}
The discretization of the $L^2$-gradient flow formulation involves several pre- and postprocessing steps. First the initial diffeomorphism $\Phi_0$ has to be computed
for a given initial density $\rho_0$. The postprocessing step corresponds to calculate the density $\rho_T$ from the final diffeomorphism $\Phi_T$ . We shall
detail these two steps in the following.

\subsection{Pre-processing: }\label{s:prepro} Let $ \rho_0(x)$ be a given smooth initial density with  $\int_{\Omega} \rho_0(x) dx = M$ and denote by $\Phi_0= \Phi_0(x) := \Phi(0,x)$ the corresponding diffeomorphism. Then the initial diffeomorphism $\Phi_0: \bar{\tilde{\Omega}} \rightarrow \bar{\Omega}$ satisfies
\begin{align}\label{e:transf}
\rho_0(\Phi_0(x)) \det D\Phi_0(x) = 1.
\end{align}
Depending on the discretization of $\tilde{\Omega}$ different approaches can be used to determine $\Phi_0(x)$ from a given $\rho_0(x)$. In the case of an equidistant mesh of squares or rectangles, one can calculate the initial diffeomorphism by solving a one-dimensional Monge-Kantorovich problem in $x_1$ direction and subsequently
a family of Monge-Kantorovich problems in the $x_2$ direction, cf. \cite{HT2001,carrillo2009numerical}. 

However this approach is not possible in case of
general quadrilateral or triangular meshes. In this case different strategies can be used: either based on the Monge Ampere
equation (giving the optimal transportation plan in case of quadratic cost), Knote theory or density equalizing maps. We shall follow the latter, which is based on an idea of Moser \cite{M1965} that
was further studied by Avinyo and co-workers in \cite{ASMV2003}. It is based on the idea of
constructing the initial diffeomorphism from solutions of the heat equation with homogeneous Neumann boundary conditions. This approach was also used in cartography to calculate density equalizing maps, cf. \cite{GN2004}. Its advantage is its flexibility - it can be used for general computational domains and their respective discretizations since
it only requires the efficient solution of the heat equation. 

Consider the heat equation written as a continuity equation on a bounded domain $\Omega \subset \mathbb{R}^2$ as
\begin{align}\label{e:cons}
\frac{\partial \rho}{\partial t} + \Div(\rho v) = 0, \qquad \mbox{with } v  = -\frac{\nabla \rho}{\rho},
\end{align}
and initial datum $\rho(0,x) = \rho_0(x)$ as well as homogeneous Neumann boundary conditions. Equation \eqref{e:cons} corresponds to the time evolution of $\rho_0$ transported by the velocity field $v  = -\frac{\nabla \rho}{\rho}$ towards the constant density $\bar{\rho} = \frac{1}{\lvert \Omega \rvert} \int_{\Omega} \rho_0(x) dx$ as $t \rightarrow \infty$.
Hence the velocity field can be calculated by solving the heat equation until equilibration.
Then the cumulative displacement $\mathbf{x}(t)$ of any point at time $t$ is determined by integrating the velocity field, which corresponds to solving
\begin{equation*}%\label{equ:displacements}
\mathbf{x}(t)=\mathbf{x}(0)+\int_0^t v(t',\mathbf{x}(t'))\,dt'.
\end{equation*}
As $t\rightarrow \infty$, the set of such displacements for all points $x = x(t)$ in $\Omega$, that is the grid points of the computational mesh, defines the new density-equalized domain. Note that we actually need to determine its inverse, since we need to find  the map which maps the constant density to the initial density $\rho_0$.  Hence we solve: 
\begin{subequations}\label{e:back}
\begin{align} 
\mathbf{x}'(t) &= v(t',\mathbf{x}(t')) = -\frac{\nabla \rho(t,\mathbf{x}(t))}{\rho(t,\mathbf{x}(t))}\\
\mathbf{x}(T) &= x_1, 
\end{align}
\end{subequations}
for all mesh points $x_1\in \Omega$. This corresponds to the solution of the integral equation $\mathbf{x}(0) = \mathbf{x}(T) - \int_0^T v(t,\mathbf{x}(t)) dt.$ Altogether the pre-processing consists of two steps:
\begin{enumerate}
\item Solve the heat equation with Neumann boundary conditions for an initial datum $\rho_0 = \rho_0(x)$ until equilibration. The time of equilibration corresponds to the time where the $L^2-$norm between two consecutive time steps is less than $10^{-6}$.
\item Starting with a given mesh at time $t=T$ calculate the initial diffeomorphism by solving \eqref{e:back} backward in time. Then the initial
diffeomorphism is given by $\Phi_0(x_1) = \mathbf{x}(0)$.
\end{enumerate}

\begin{rem}
Note that the one could also use the initial density $\rho_0$ as a reference measure, then the initial diffeomorphism $\phi_0$ is the identity map. In this case equation \eqref{e:varchange} would read as 
\begin{align*}
\rho(\Phi(x)) \det(D \Phi(x)) = \rho_0(x).
\end{align*}
\end{rem}

\subsection{Post-processing:} The post-processing step corresponds to calculate $\rho_T := \rho(t = T,x)$ given the final diffeomorphism $\Phi_T := \Phi(t = T,x)$.
Since we expect numerical artifacts  in case of compactly supported and measure valued solutions we solve a 
regularized version of \eqref{e:varchange} (in 2D) given by:
\begin{align}\label{e:regvarchange}
\varepsilon \Delta \rho_T(\Phi_T(x)) + \rho_T(\Phi_T(x)) = \frac{1}{\det D\Phi_T(x)} \text{ with } 0 < \varepsilon \ll 1.
\end{align}
Note that equation \eqref{e:regvarchange} can be easily implemented using finite elements, since the assembling of the system matrix is based on the
transformation of the triangle to the reference triangle. In case of equation \eqref{e:regvarchange} we have to perform two
 successive transformations. First the transformation to the displaced element given by $\Phi$ and second to the reference element.
 
%%%%%%%%%%%%%%%%%%%%%%%%%%%%%%%%%%%%

\section{Numerical simulations}\label{s:numerics}

In this section we illustrate the behavior of the numerical solver with simulations in spatial dimension one and two. The 1D simulations are based on a finite difference discretization, in 2D we use finite elements. The respective solvers are based on Matlab in 1D and on the finite element package Netgen/NgSolve in 2D. In the following we illustrate the flexibility of our approach with various simulations for a large class of PDEs. We start by summarizing all steps of the numerical solver in Algorithm \ref{alg}.
\begin{algorithm}
\caption{Let  $\rho_0 = \rho_0(x)$ be a given initial density with $\int_{\Omega} \rho_0(x) dx = M$ and $\Omega = [a,b]$ or $\Omega = [a,b]\times[c,d]$.} 
\label{alg}
\begin{enumerate}
\item Determine the initial diffeomorphism 
\begin{enumerate}
\item 1D: Define the initial diffeomorphism $\Phi_0: \tilde{\Omega}=[0,M]\to [a,b]$ by: for every discretization point $x_i \in [0,M]$ solve the one-dimensional Monge-Kantorovich problem for $\Phi_i = \Phi(x_i)$:
\begin{align*}
\int_{a}^{\Phi_i} \rho_0(y) dy = x_i,
\end{align*}
using Newton's method.
\item 2D: Follow the two-step algorithm outlined in Section \ref{s:prepro}:
Solve the heat equation $\partial_t \rho = \Delta \rho$ with homogeneous Neumann boundary conditions using an $H^1$ conforming finite element method
of order $p=6$ until equilibration at time $t=T$. Determine the initial diffeomorphism by transporting each mesh point of the computational mesh $V(T) = (\mathbf{x}_1^V(T), \mathbf{x}_2^V(T)) \in \mathcal{V}$ backward in time
\begin{align*}
\mathbf{x}(0) = \mathbf{x}(T) - \int_0^T v(t,\mathbf{x}(t)) dt, ~~v(t,\mathbf{x}(t)) = -\frac{\nabla \rho(t,\mathbf{x}(t))}{\rho(t,\mathbf{x}(t))}.
\end{align*}
\end{enumerate}
\item At every time step $t^{n+1} = (n+1) \Delta t$ in $t \in (0,T]$ solve 
\begin{enumerate}
\item 1D: the implicit equation \eqref{e:implicitdiffeo} using Newton's method and finite difference discretization.
\item 2D: the nonlinear operator equation $F(\Phi, \varphi) = 0$, given by \eqref{e:F} using Newton's method and a spatial discretization of $H^1$ conforming finite elements of order $p=1$.
\end{enumerate}
\item Recover the final density $\rho= \rho(T)$ by 
\begin{enumerate}
\item 1D: calculating the final density at every discrete point $\Phi_i = \Phi(x_i)$ from the relation
\begin{align*}
\rho(\Phi(x)) = \frac{1}{\det D \Phi(x)}.
\end{align*}
\item 2D: solving the regularized equation
\begin{align*}
\varepsilon \Delta \rho(\Phi(x)) + \rho(\Phi(x)) = \frac{1}{\det D \Phi(x)}
\end{align*}
on the transformed mesh $\mathcal{N}$ using $H^1$ conforming elements of order $p=4$.
\end{enumerate}
\end{enumerate}
\end{algorithm}

\subsection{Numerical simulations in spatial dimension one}\label{num_sim1D}

If not stated otherwise we initiate all numerical simulations with a Gaussian of the form:
\begin{align}\label{e:gaussian}
\rho_0(x) =\frac{M}{\sqrt{2\pi\sigma^2}}e^{\frac{-x^2}{2\sigma^2}}
\end{align}
with $M, \sigma >0$. Note that $\int_{-\infty}^{\infty} \rho_0(x) dx = M$. We perform our simulations on a bounded domain $\tilde{\Omega}=[0,\tilde{M}]$, where $\tilde{M}$ is the approximate value of $\int_{\Omega} \rho_0(x) dx$.
Hence, the initial diffeomorphism $\Phi_0$ is defined on $\tilde{\Omega}$.

\subsubsection{Nonlinear diffusion equations} 
\begin{figure}[h!]
\begin{center}
\subfloat[Reconstructed BP profile  $\rho = \rho(x,t)$ at time $t=0.021$.]{\includegraphics[width=0.42\textwidth]{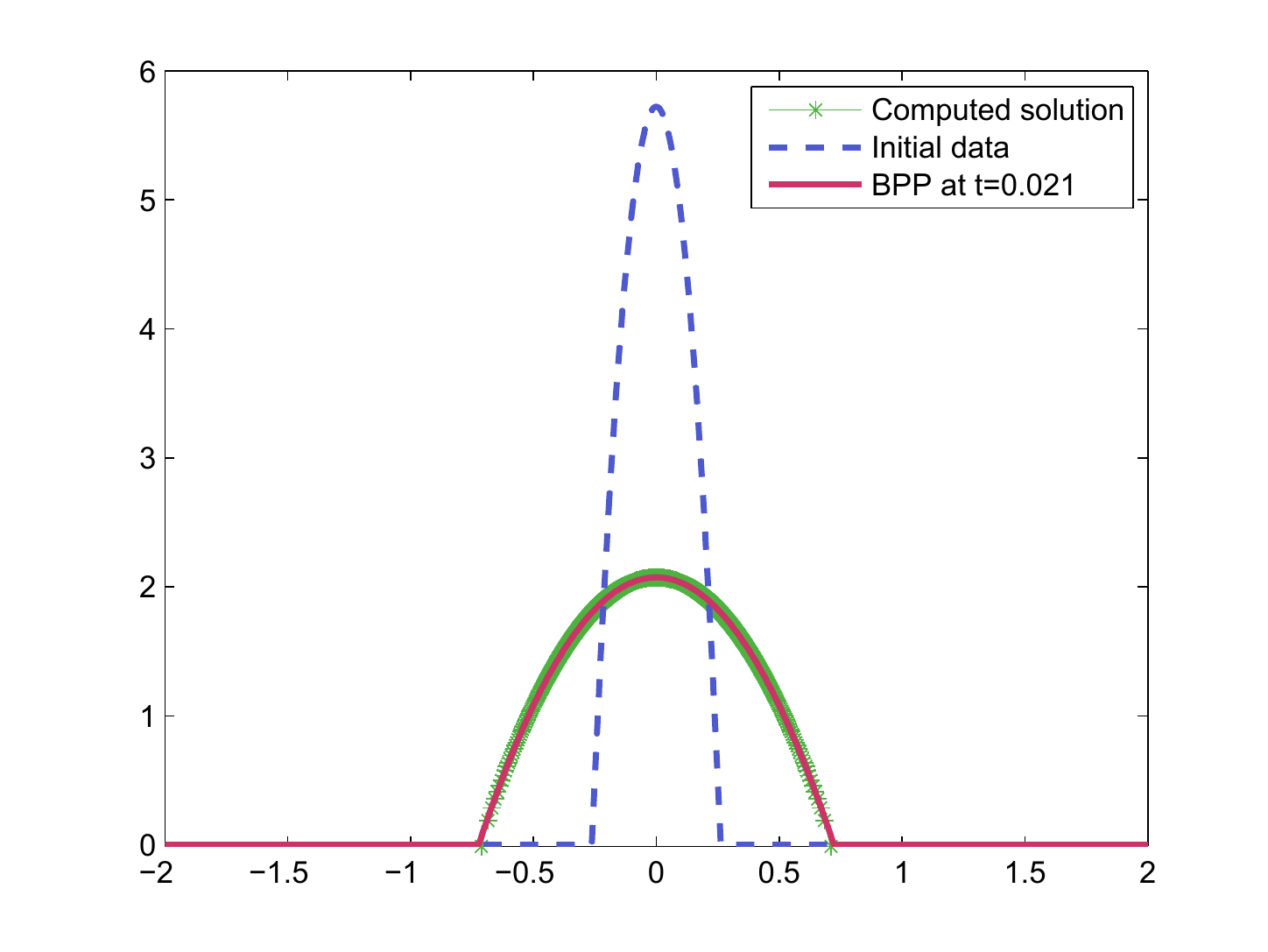}}\hspace*{0.5cm}
\subfloat[Evolution of the entropy and the expected decay of $t^{-\alpha}$ in time with a logarithmic scale.]{\label{f:pme1_entropy}\includegraphics[width=0.42\textwidth]{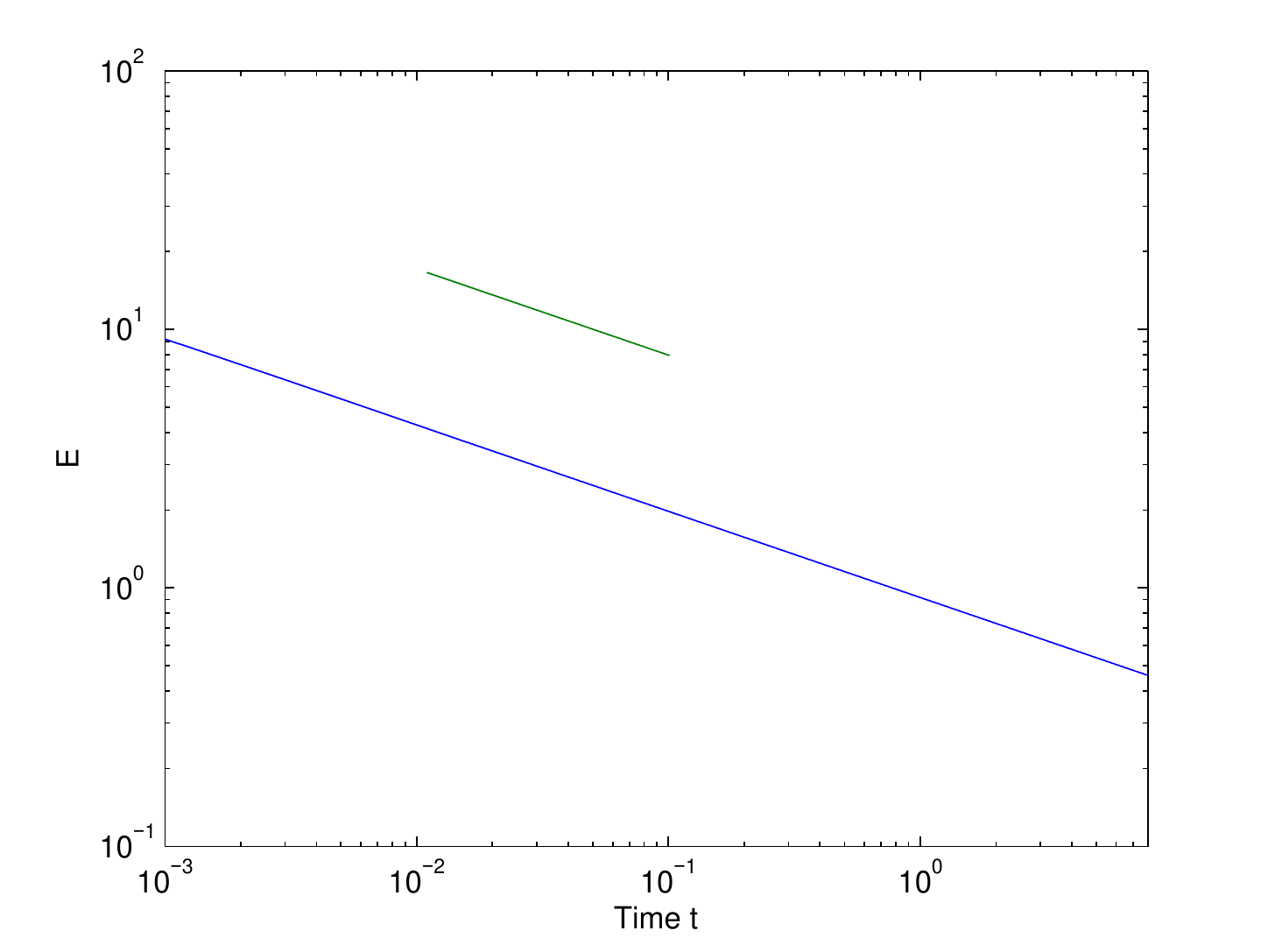}}
\caption{Solution of the PME for $m=2$ after 2000 time steps.} \label{f:pme1}
\end{center}
\end{figure}

\begin{figure}[h!]
\begin{center}
  \subfloat[Reconstructed BP profile $\rho= \rho(x,t)$ at time $t = 0.021$.]{\includegraphics[width=0.42\textwidth]{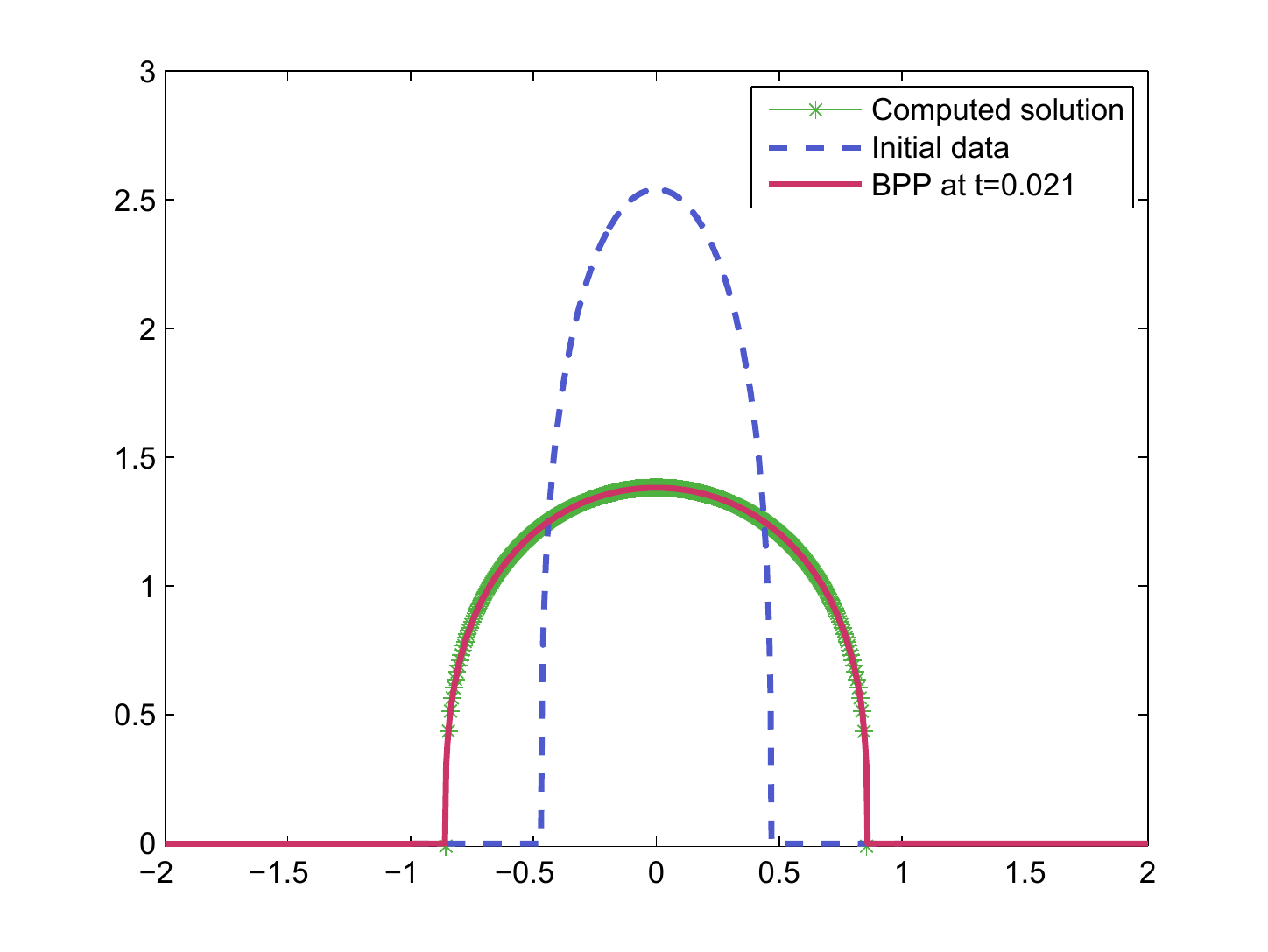}}\hspace*{0.5cm}
\subfloat[Evolution of the entropy and the expected decay of $t^{-3\alpha}$ in time with a logarithmic scale.]{\label{f:pme2_entropy}\includegraphics[width=0.42\textwidth]{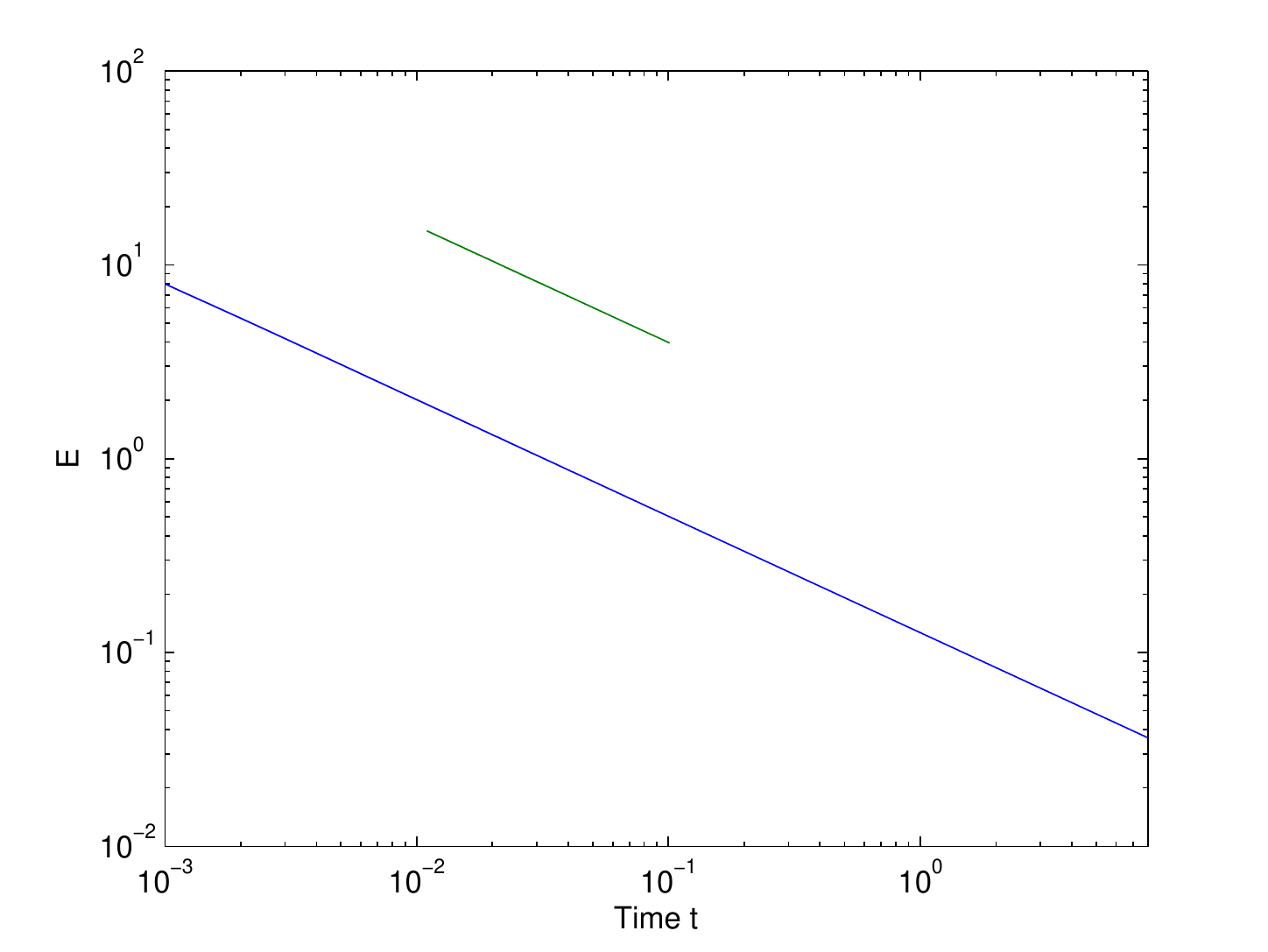}}
\caption{Solution of the PME for $m=4$ after 2000 time steps.} \label{f:pme2}
\end{center}
\end{figure}

\noindent In our first example we illustrate the behavior of our scheme for the Porous medium equation, that is $U(s) = \frac{1}{m-1} s^m, m > 1$ and $V=W=0$. Let $\Omega = [-1,1]$ denote the computational domain discretized into $501$ intervals of size $\Delta x = \frac{2}{501}$. The initial density is given by a Barenblatt Pattle profile (BPP) at time $t_0=10^{-3}$, i.e.
\begin{align*}
\rho_0(x) = \frac{1}{t_0^\alpha}\left(c-\alpha\frac{m-1}{2m}\frac{x^2}{t_0^{2\alpha}}\right)_+^{\frac{1}{m-1}},
\end{align*}
where $\alpha=\frac{1}{m+1}$ and $c$ is chosen such that $\int_{-1}^1 \rho_0(x)\,dx=2$.
The discrete time steps are set to $\Delta t=10^{-5}$.
We consider the cases $m=2$ and $m=4$. Figures \ref{f:pme1} and \ref{f:pme2} show the density profiles $\rho$ at time $t=0.021$ as well as the evolution of the free energy in time, which decays like $t^{-\alpha(m-1)}$, cf. \cite{vazquez2006smoothing}. As seen from Figures \ref{f:pme1_entropy} and \ref{f:pme2_entropy}, these decays (indicated by the green lines) are perfectly captured by the scheme validating the chosen discretization.  Note that the boundary points of the support are calculated using an explicit one-sided difference scheme of the velocity as proposed by Budd et al., see \cite{budd1999self}. In particular, writing the porous medium equation as 
\[\rho_t+\nabla \cdot (\rho v)=0\] 
with the velocity $v=-\frac{m}{m-1}(\rho^{m-1})_x$ leads to the approximation $\Phi^{n+1}=\Phi^n+\Delta t v$ at the boundary points. When determining $\rho_x$ we assume that the density at the boundary of the support equals zero.
 

\subsubsection{Fokker-Planck type equations}
\noindent Next we study the behavior of a nonlinear Fokker-Planck equation in the case of an internal energy $U(s)=\nu \frac{ s^m}{m}$, $m>1$ and a double-well confining potential $V(x)=\frac{x^4}{4}-\frac{x^2}{2}$. Let $\Omega = [-2,2]$ be the computational domain discretized into $501$ intervals of size $\Delta x = \frac{4}{501}$. The initial density is given by \eqref{e:gaussian} with $\sigma=0.2$ and $M=1$. 
\begin{figure}[h!]
\begin{center}
\subfloat[Density $\rho$ at time $T=9.718$.]{\includegraphics[width=0.42\textwidth]{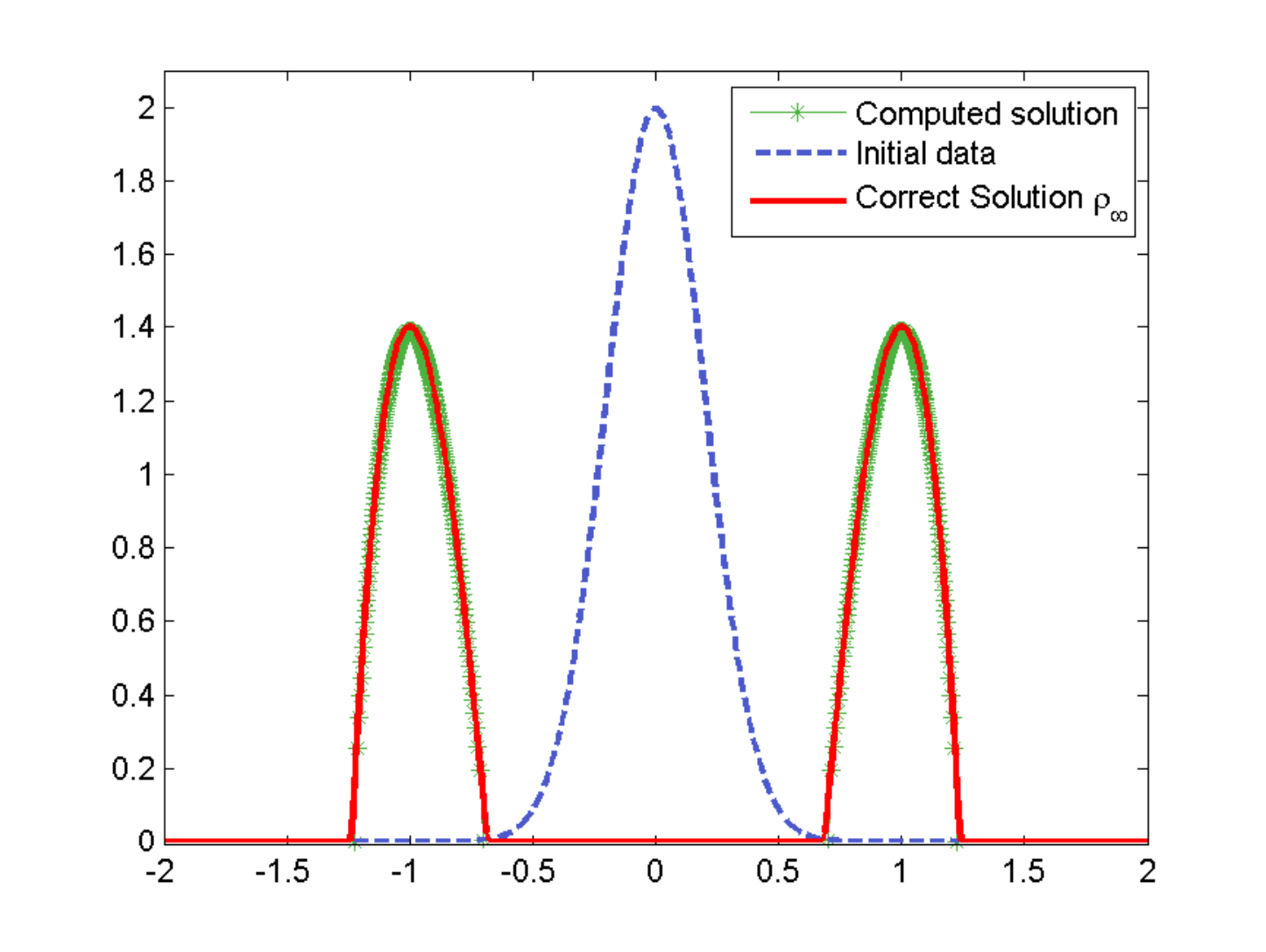}} \hspace*{0.5cm}
\subfloat[Evolution of the entropy in time with a logarithmic scale for the y-axis.]{\label{f:ex7_entropy}\includegraphics[width=0.42\textwidth]{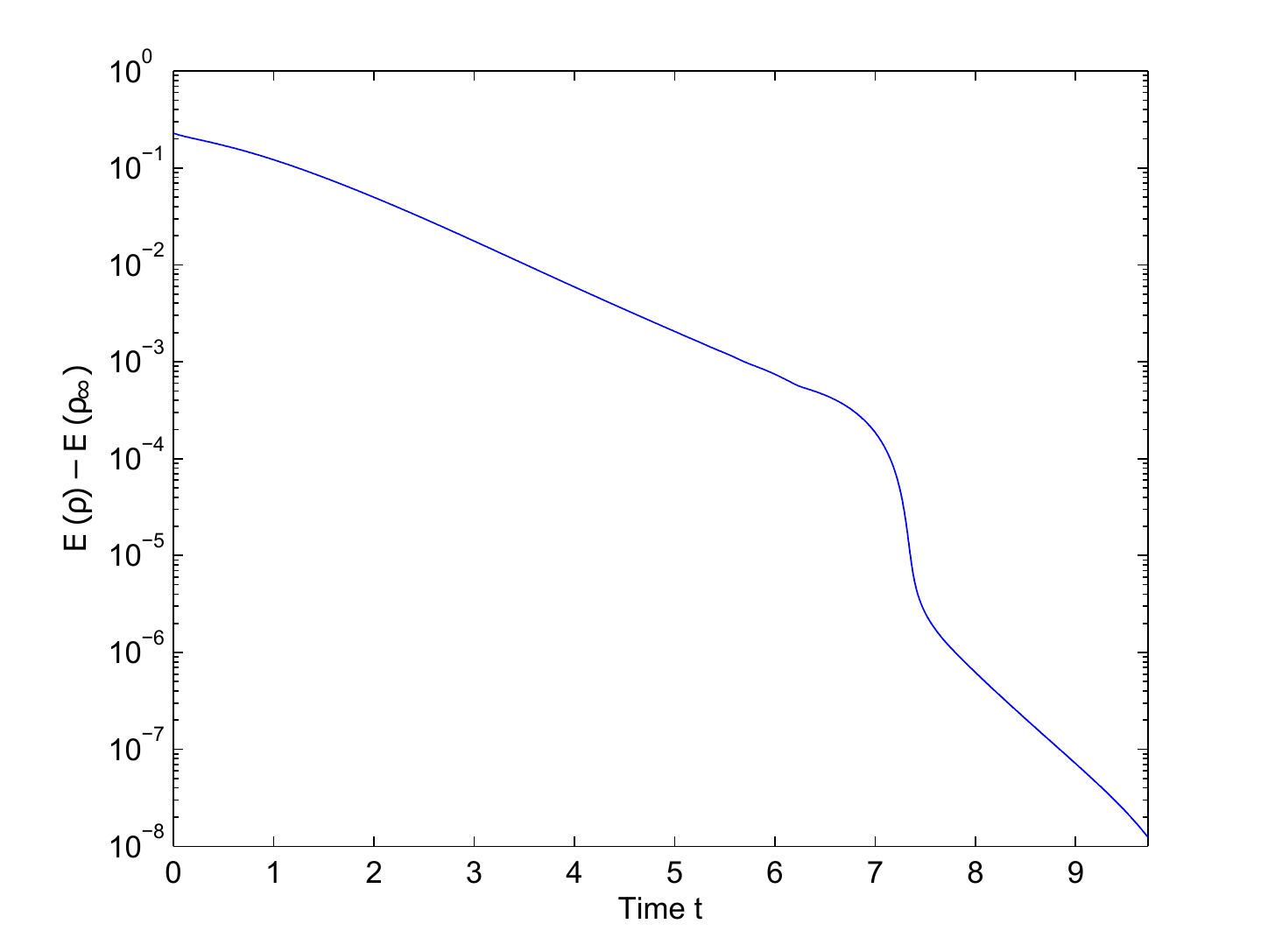}}
\caption{Solution of the nonlinear Fokker-Planck equation for $m=2$ and $\nu=0.05$.} \label{f:ex7}
\end{center}
\end{figure}
The discrete time steps are set to $\Delta t=10^{-3}$. In this case steady states are given by
\begin{align}\label{steady_state}
\rho_\infty(x)=\frac1{\nu} \left(C(x)-V(x)\right)_+^{\frac{1}{m-1}},
\end{align}
where $C = C(x)$ is a piecewise constant function possibly taking different values on each connected part of the support, cf. \cite{CCH2015}. The density at time $T$ and the entropy decay with respect to $\rho_\infty$ are illustrated in Figure \ref{f:ex7}. We always choose $T$ in Section \ref{num_sim1D} as the time where the $L^2-$norm of the difference of two consecutive diffeomorphisms is less than $10^{-6}$ and refer to $\rho$ at time $T+t_0$ for $t_0=500\Delta t$ as the steady state $\rho_\infty$. We observe that the stationary state has two connected components. As seen in Figure \ref{f:ex7_entropy}, the entropy decays abruptly once the support separates into two pieces before final convergence towards the steady state with the lowest free energy filling with equal mass each of the two wells of the potential, i.e., $C=C(x)$ is equal in each connected component.

\begin{figure}[h!]
\begin{center}
\subfloat[Density $\rho$ at time $T=3.172$.]{\label{f:ex7_s}\includegraphics[width=0.42\textwidth]{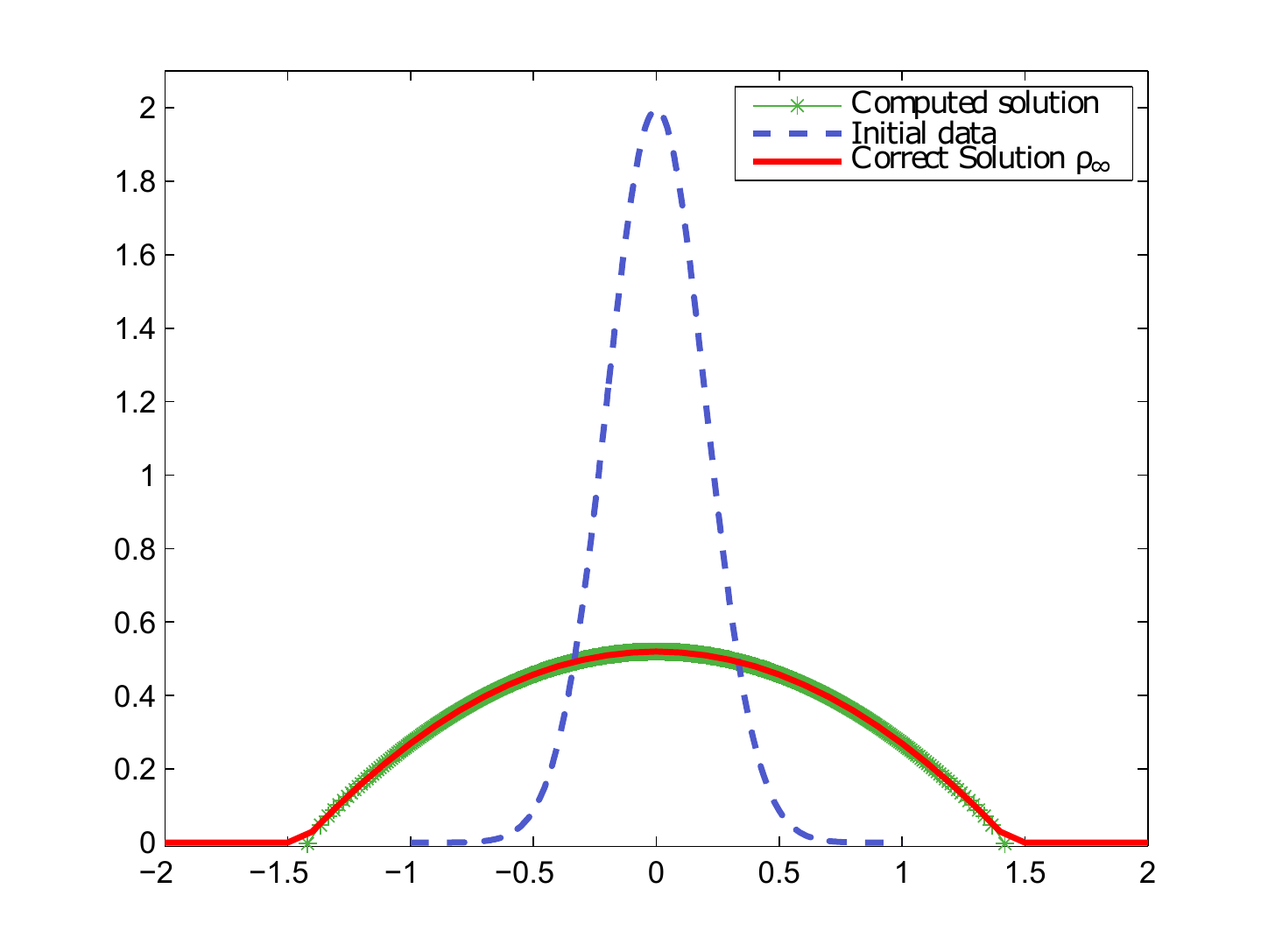}} \hspace*{0.5cm}
\subfloat[Semi-log error plot at time t for different values for $N=\frac{|\Omega|}{\Delta x}$ with respect to the steady state.]{\label{f:ex7_e} \includegraphics[width=0.42\textwidth]{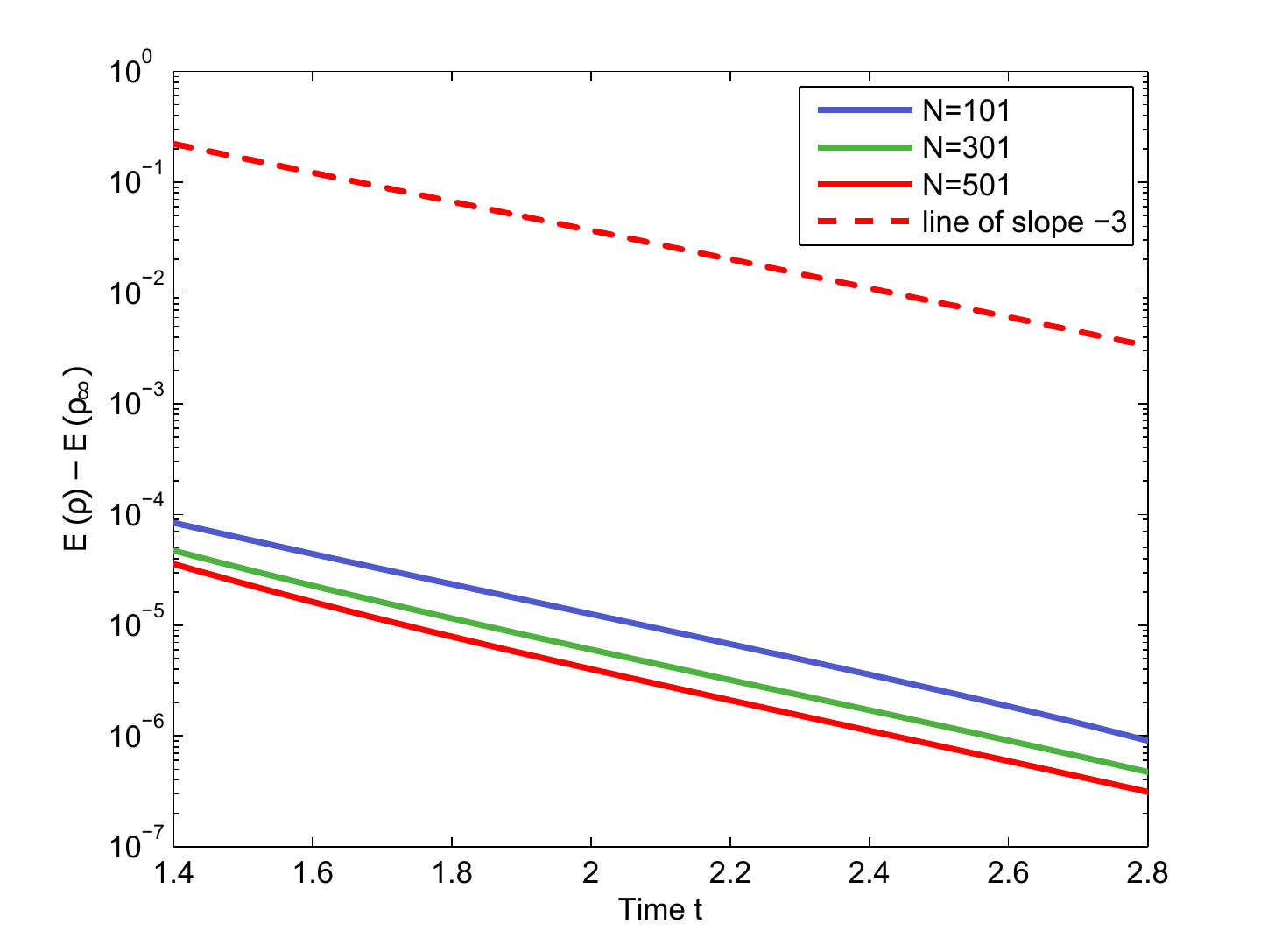}} 
\caption{Solution of the nonlinear Fokker-Planck equation for $m=2$ and $\nu=2$.}
\end{center}
\end{figure}

If the confining potential is chosen as the harmonic potential $V(x)=\frac{x^2}{2}$, the steady state is also given by \eqref{steady_state}, cf. \cite{carrillo2007strict}. The convergence towards the steady state is exponential, more precisely Carrillo et al.  \cite{carrillo2007strict} showed that the distance towards equilibrium, i.e. $d_W^2(\rho(t),\rho_\infty)$, converges like $\mathcal{O}(e^{-(m+1)t})$ since the initial data has zero center of mass. The same convergence behavior can be observed in the numerical simulations, see Figure \ref{f:ex7_s} and \ref{f:ex7_e}. Note that the rate of convergence in relative energy in Figure \ref{f:ex7_e} coincides with the rate of convergence of $d_W^2(\rho(t),\rho_\infty)$.

\subsubsection{Aggregation equations}
\noindent Next we consider aggregation potentials of the form 
\begin{align}\label{e:W2}
W(x)=\frac{|x|^a}{a}-\frac{|x|^b}{b},
\end{align}
for  $a>b\geq 0$ using the convention that $|x|^0/0=\ln|x|$. The set of stationary states of these equations can be quite complex, see \cite{Bertozzietal,BCLR,BCLR2}, depending on the dimension. In one dimension with $a=2$ and $b=1$ the steady state profiles are constant on an interval, cf. \cite{fellner2010stable,fellner2011stability}, for $a=2$, $b=0$ they correspond to the semicircle law, cf. \cite{opac-b1093628,carrillo2012mass}.

\begin{figure}[h!]
\begin{center}
\subfloat[Density $\rho$ at time $T = 2.412$.]{\includegraphics[width=0.42\textwidth]{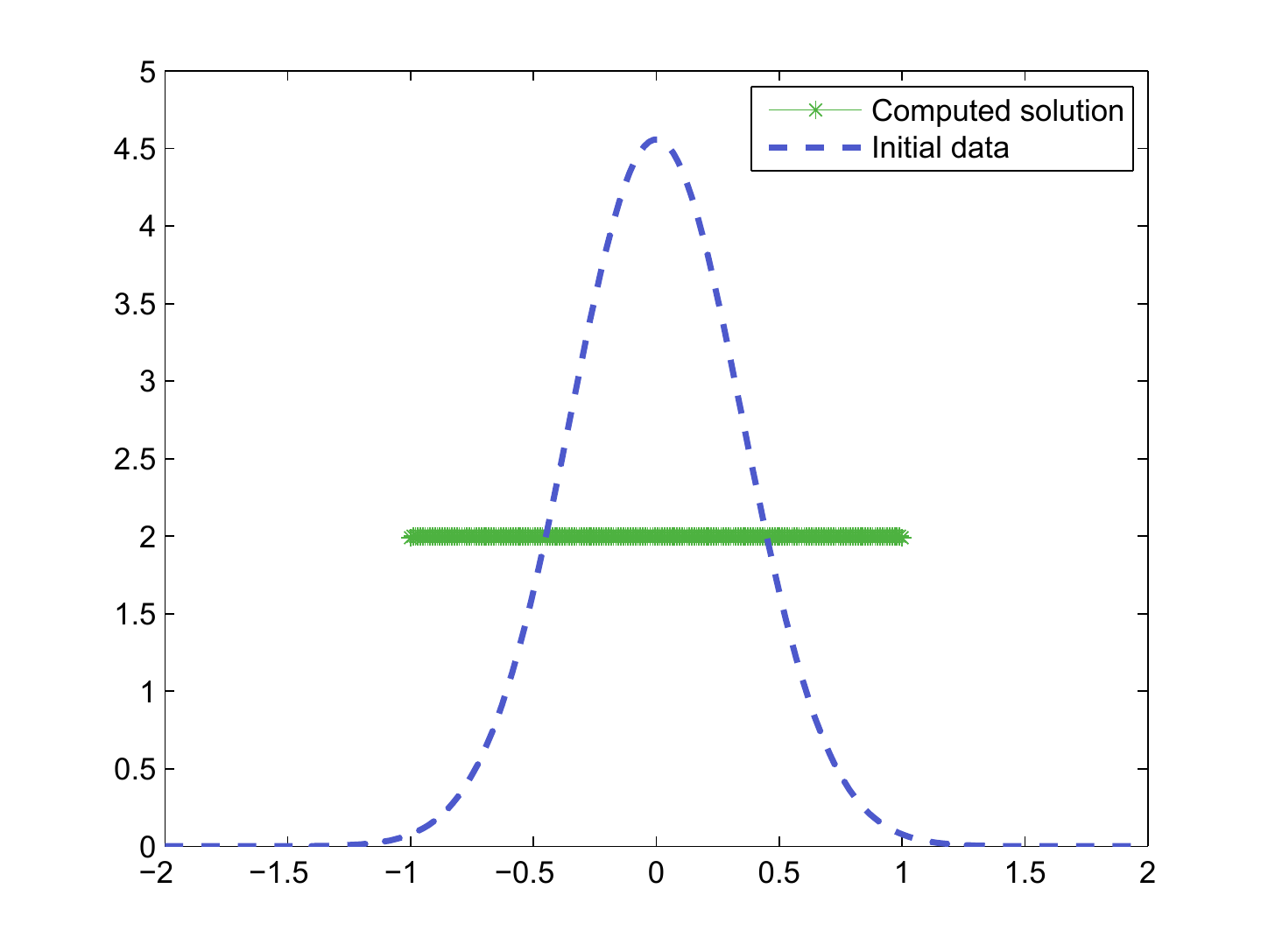}} \hspace*{0.5cm}
\subfloat[Evolution of the entropy and the decay of $e^{-7.5t}$ in time with a logarithmic scale for the y-axis.]{ \label{f:ex8_entropy}\includegraphics[width=0.42\textwidth]{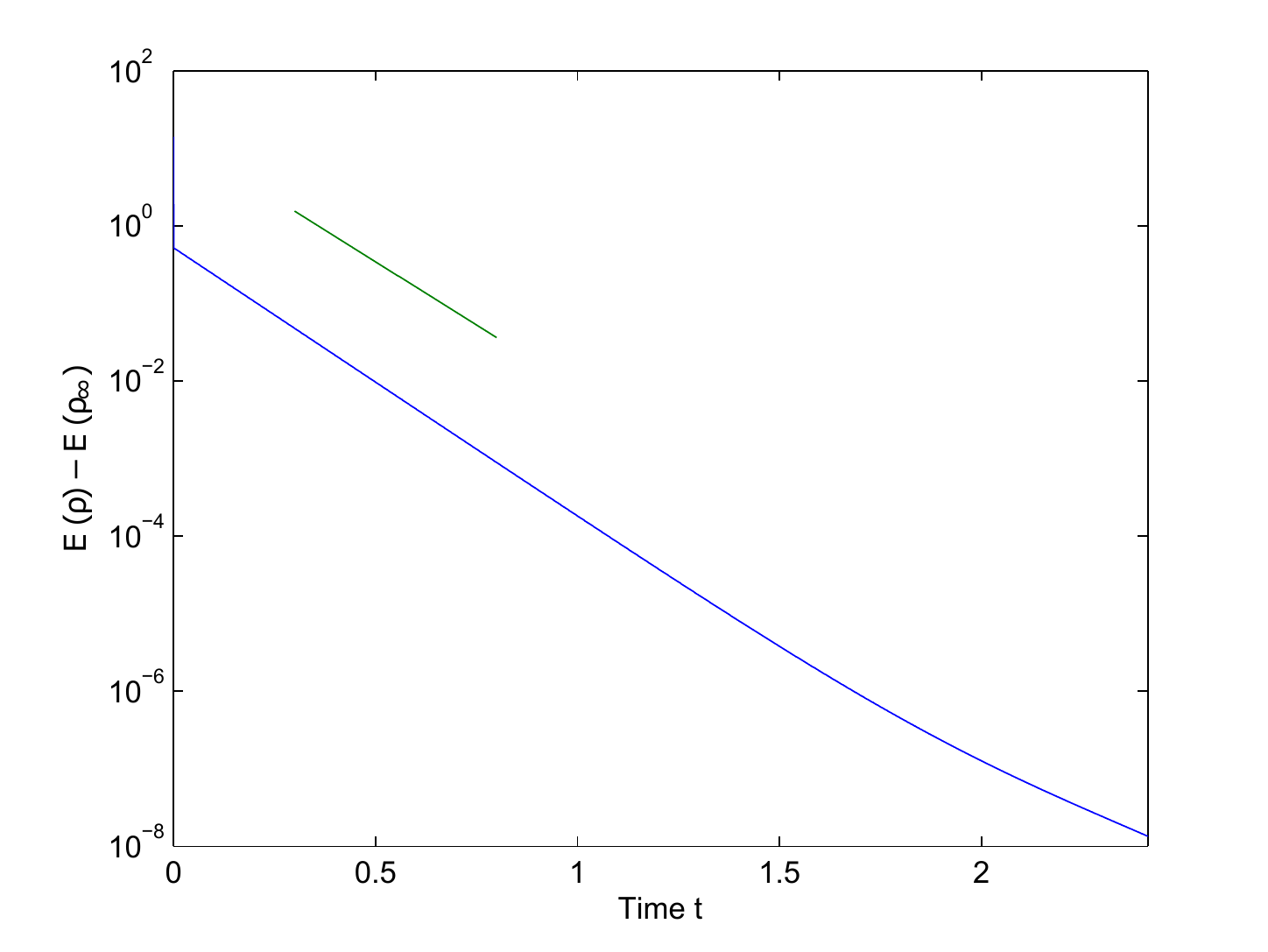}}
\caption{Solution of the aggregation equation with potential \eqref{e:W2} ($a=2$, $b=1$).} \label{f:ex8}
\end{center}
\end{figure}

\noindent Let $a=2, b=1$ and set $\Omega = [-2,2]$ split into $201$ intervals of size $\Delta x = \frac{4}{201}$. The initial datum is given by \eqref{e:gaussian} with $\sigma=0.35$ and $M=4$. The discrete time steps are set to $\Delta t=10^{-3}$. Figure \ref{f:ex8} shows the computed solution at time $T$ and the entropy decay with respect to $\rho_\infty$. The numerical simulations confirm the theoretical results, i.e. the computed stationary profile corresponds to the constant density $\rho_{\infty}(x) = 2$ for all $x \in [-1,1]$. Note that we use the mid-point rule to calculate the convolution integral in \eqref{e:implicitdiffeo} on the boundary elements, as proposed in  \cite{CCH2015}. 

\begin{figure}[h!]
\begin{center}
\subfloat[Density $\rho$ at time $T = 4.891$.]{\includegraphics[width=0.42\textwidth]{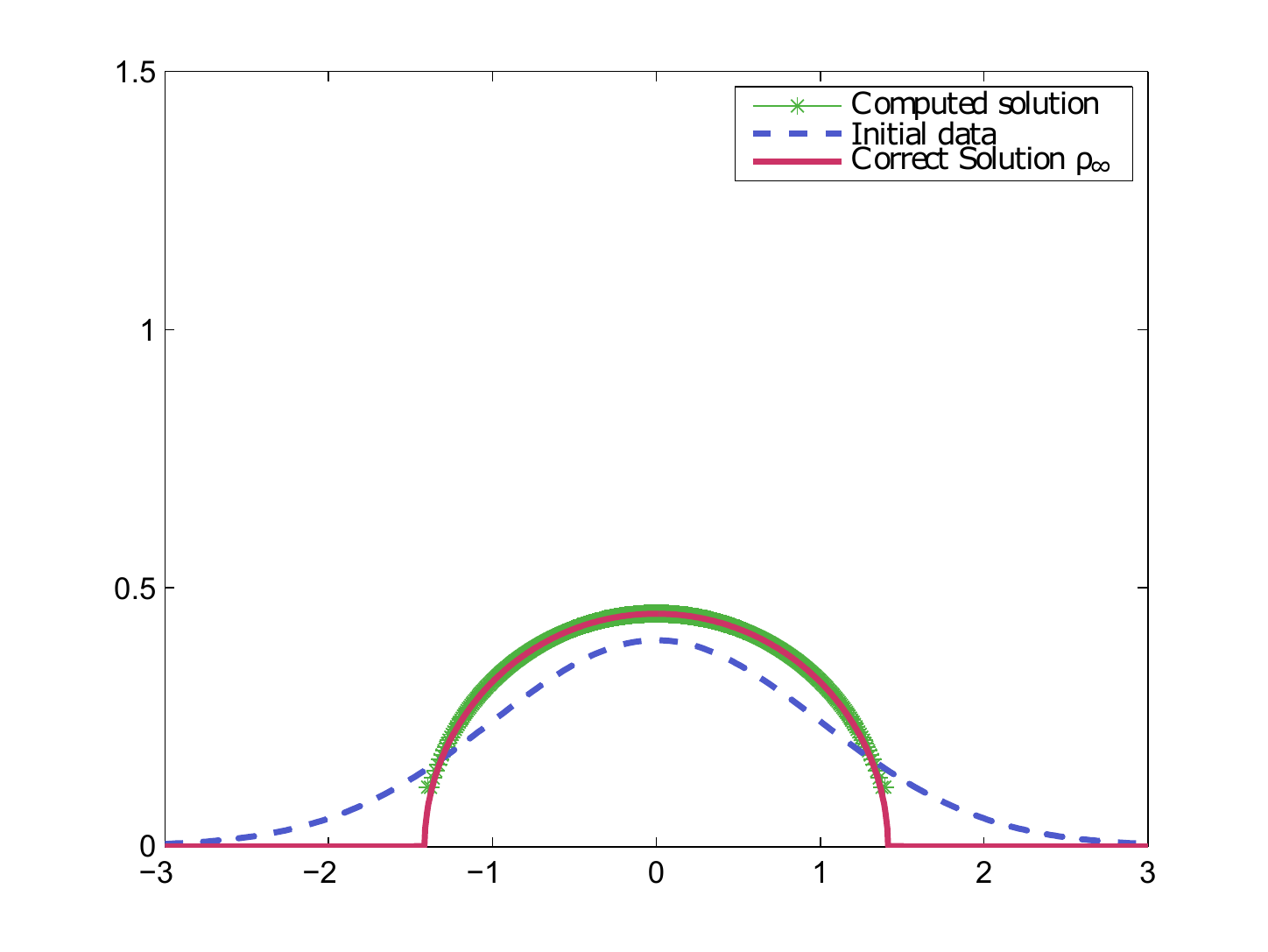}} \hspace*{0.5cm}
\subfloat[Evolution of the entropy and the decay of $e^{-3t}$ in time with a logarithmic scale for the y-axis.]{\label{f:ex9_entropy}\includegraphics[width=0.42\textwidth]{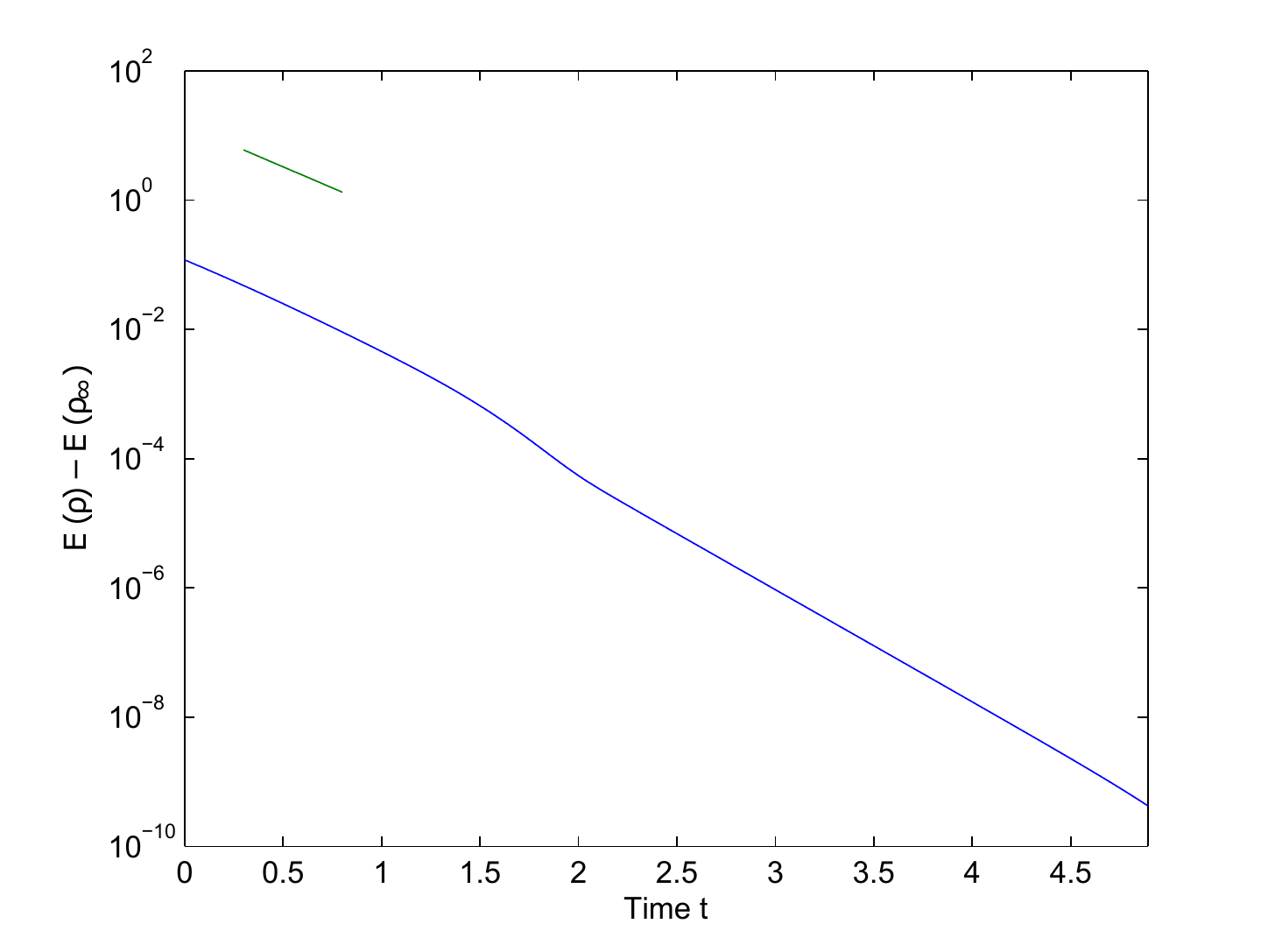}}
\caption{Solution of the aggregation equation with potential \eqref{e:W2} ($a=2$, $b=0$).} \label{f:ex9}
\end{center}
\end{figure}

In the case $a=2, b=0$ we set $\Omega=[-6,6]$ split into $501$ intervals of size $\Delta x = \frac{12}{501}$. Then the unit-mass steady state is given by 
\begin{align*}%\label{ex9_steady}
\rho_\infty(x)= \begin{cases} \frac{1}{\pi} \sqrt{2-x^2},\quad &|x|\leq \sqrt{2}, \\
0\  & \text{otherwise}. \end{cases}
\end{align*}
This result is confirmed by the numerical simulations illustrated in Figure \ref{f:ex9} for an initial datum $\rho_0$ given by \eqref{e:gaussian} with $\sigma =1$ and $M = 1$. 
It is known that the solutions decay in $d_W$ toward the corresponding stationary states $\rho_\infty$ exponentially fast. We observe that behavior at the level of the relative energy with numerical rates of decay given by 7.5 and 3 approximately in Figures \ref{f:ex8_entropy} and  \ref{f:ex9_entropy}.

\subsubsection{Modified Keller-Segel model}

\noindent In our final example we consider a modified Keller-Segel (KS) model. It is well known that solutions of the classical Keller-Segel may blow up in finite
time in space dimensions $d \geq 2$, see cf. \cite{MR2391528,blanchet2008convergence}. The same behavior can be observed in 1D by using the corresponding 2D interaction potential instead. Then the modified KS model reads as:
\begin{align}
\begin{aligned}
\partial_t \rho  &= \nabla \cdot (\rho (\nabla (\ln \rho+ W*\rho))),\\
W&=\frac{\chi}{d\pi}\ln |x|,\\
\rho(x,0)&=\rho_0\geq 0.
\end{aligned}
\label{e:mksm}
\end{align}
The blow-up behavior depends on the initial mass $M_0 = \int_{\Omega} \rho_0(x) dx$. If $M_0 < M_c$, where $\chi M_C = 2d^2 \pi$ denotes the critical mass, the system has a global in time solution. If $M_0 > M_c$ solutions blow up in finite time.

\begin{figure}[h!]
\begin{center}
\includegraphics[width=0.42\textwidth]{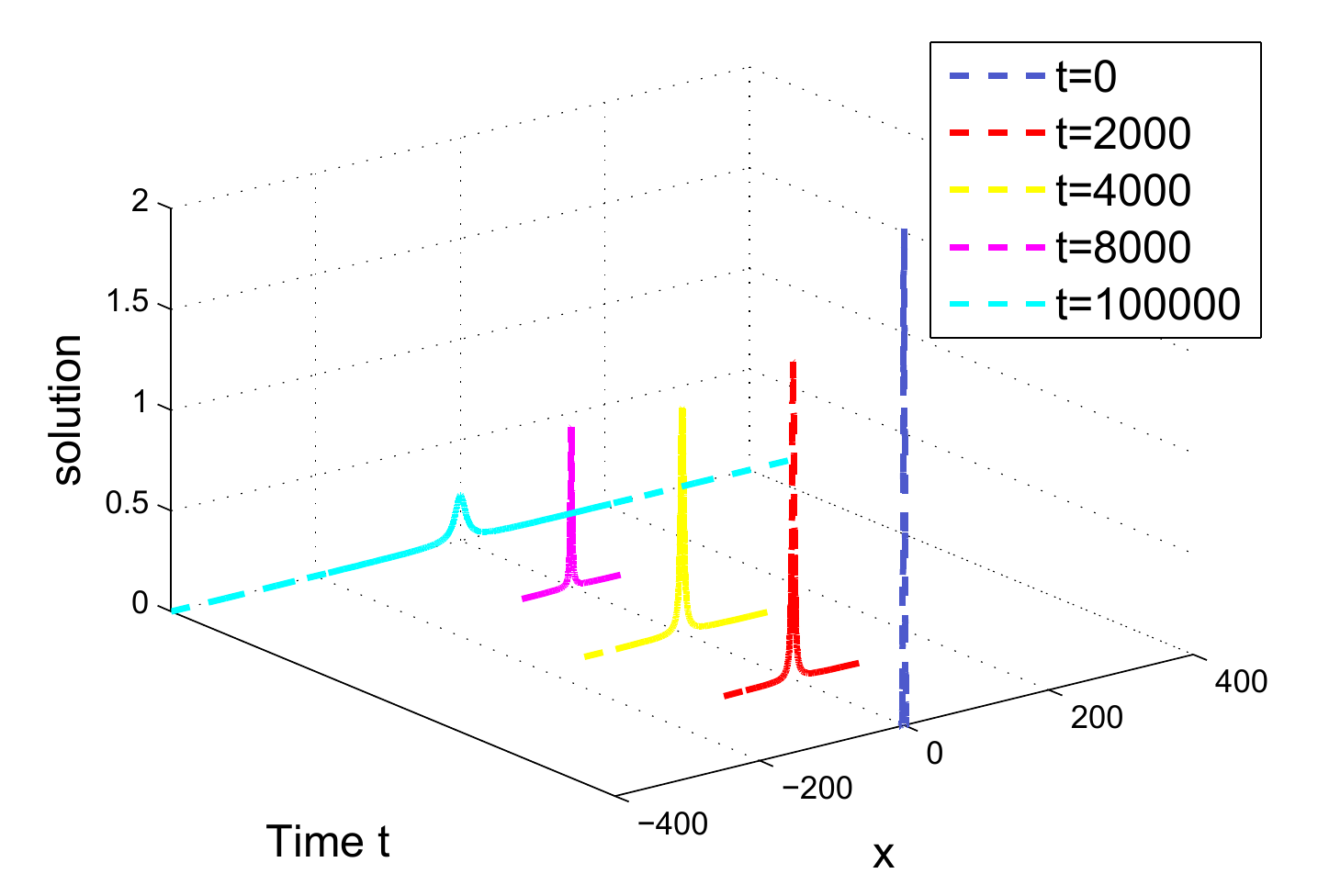} \hspace*{0.5cm}
\caption{$M=2\pi-0.1<M_c$. Evolution of the density $\rho$ for $\chi=1$.}\label{f:ex11}
\end{center}
\end{figure}

\begin{figure}[h!]
\begin{center}
\subfloat{\includegraphics[width=0.42\textwidth]{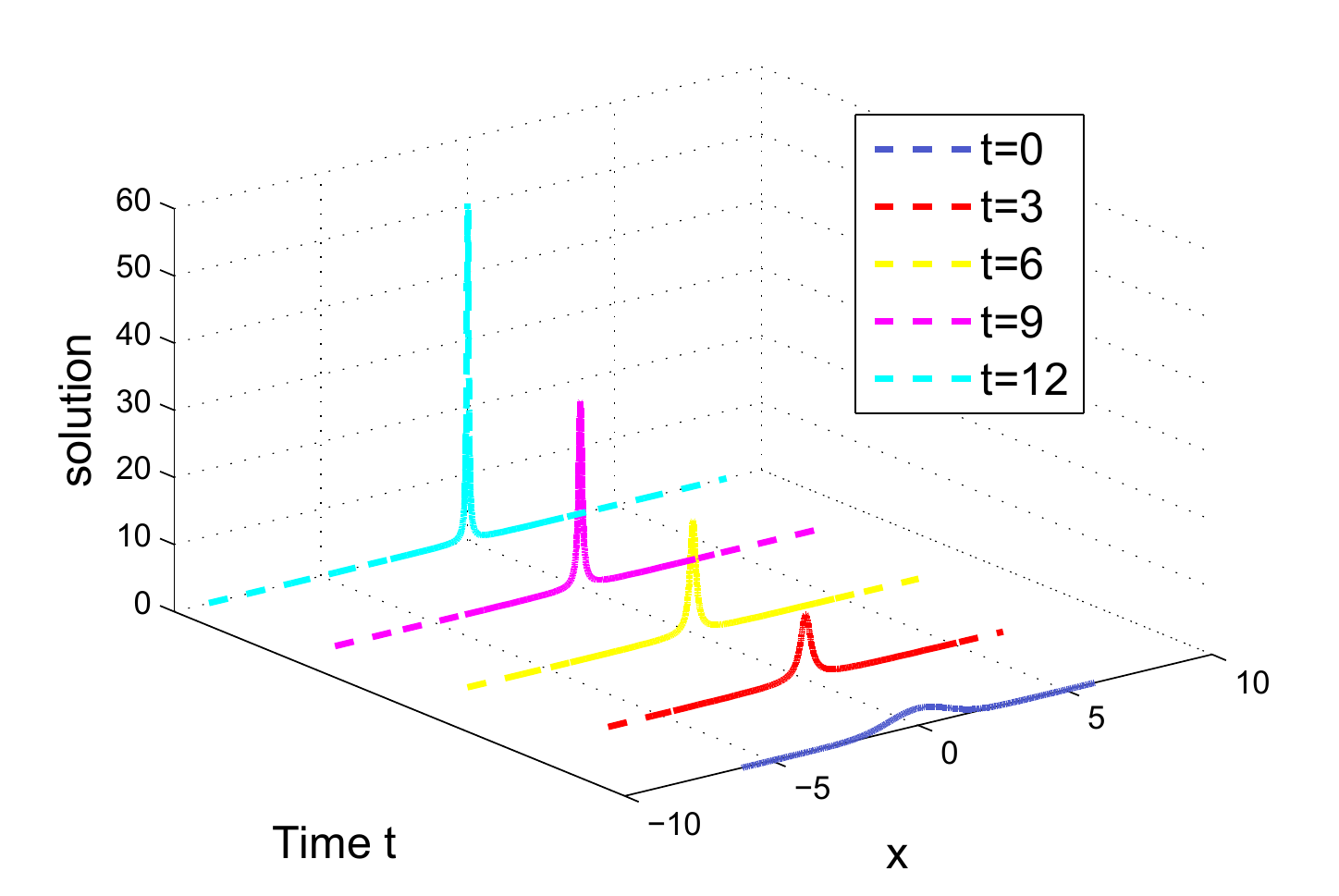}} \hspace*{0.5cm}
\subfloat{\includegraphics[width=0.42\textwidth]{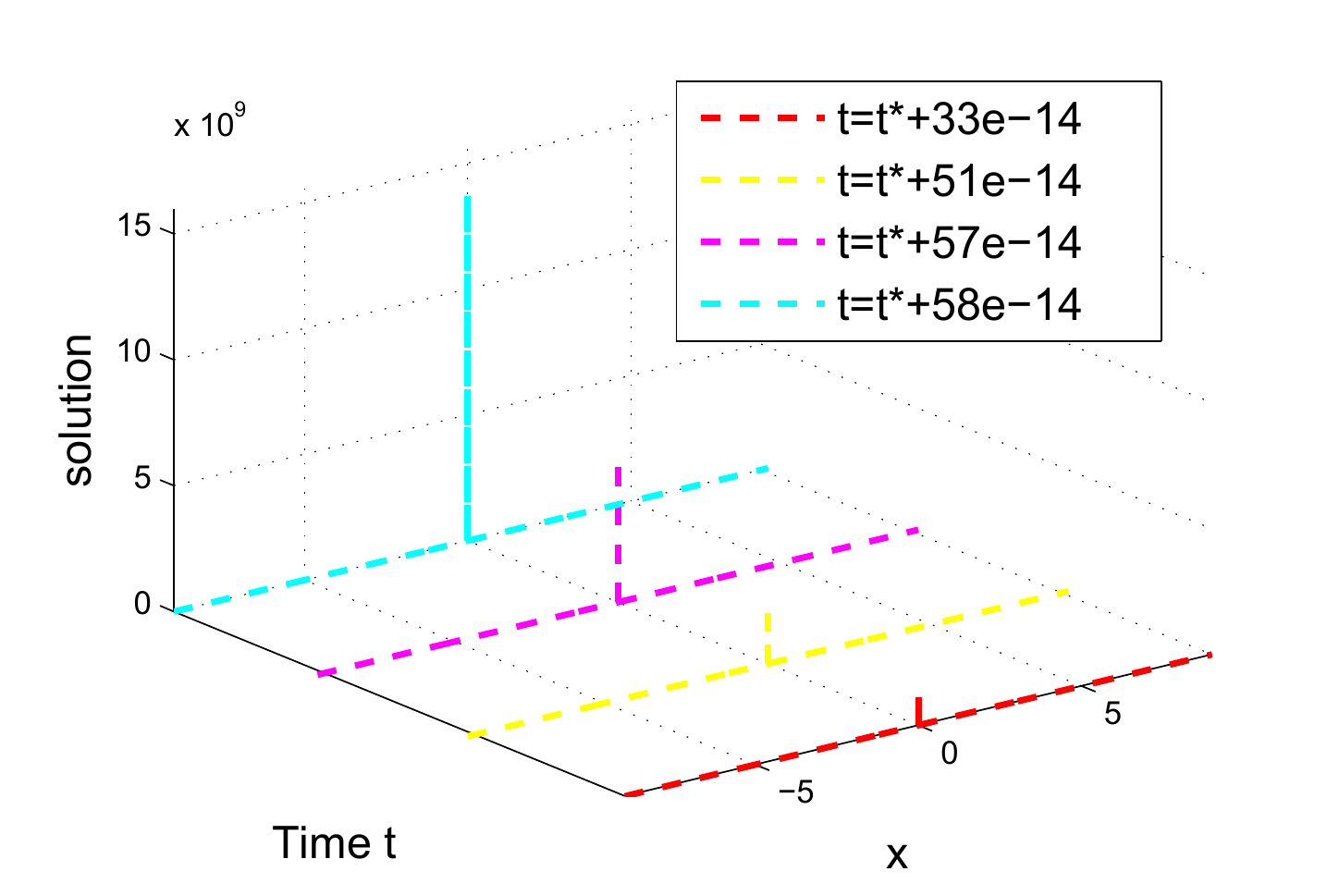}}
\caption{$M=2\pi+0.1>M_c$. Evolution of the density $\rho$ for $\chi=1$, where $t^*=16.818399684193$.}\label{f:ex10}
\end{center}
\end{figure}

Let $\Omega = [-6,6]$ be the computational domain discretized into $501$ intervals of size $\Delta x = \frac{12}{501}$. The initial density is given by \eqref{e:gaussian} with $\sigma=1$ and the discrete time steps are set to $\Delta t=10^{-1}$. If the initial data has a smaller mass, in our case $M_0 = 2\pi - 0.1$ the solution diffuses to zero, see Figure \ref{f:ex11}. If the $M_0 > M_C$ the solution blows up as illustrated in Figure \ref{f:ex10} in the case of $M_0 = 2\pi + 0.1$. 

\begin{figure}[h!]
\begin{center}
\subfloat{\includegraphics[width=0.42\textwidth]{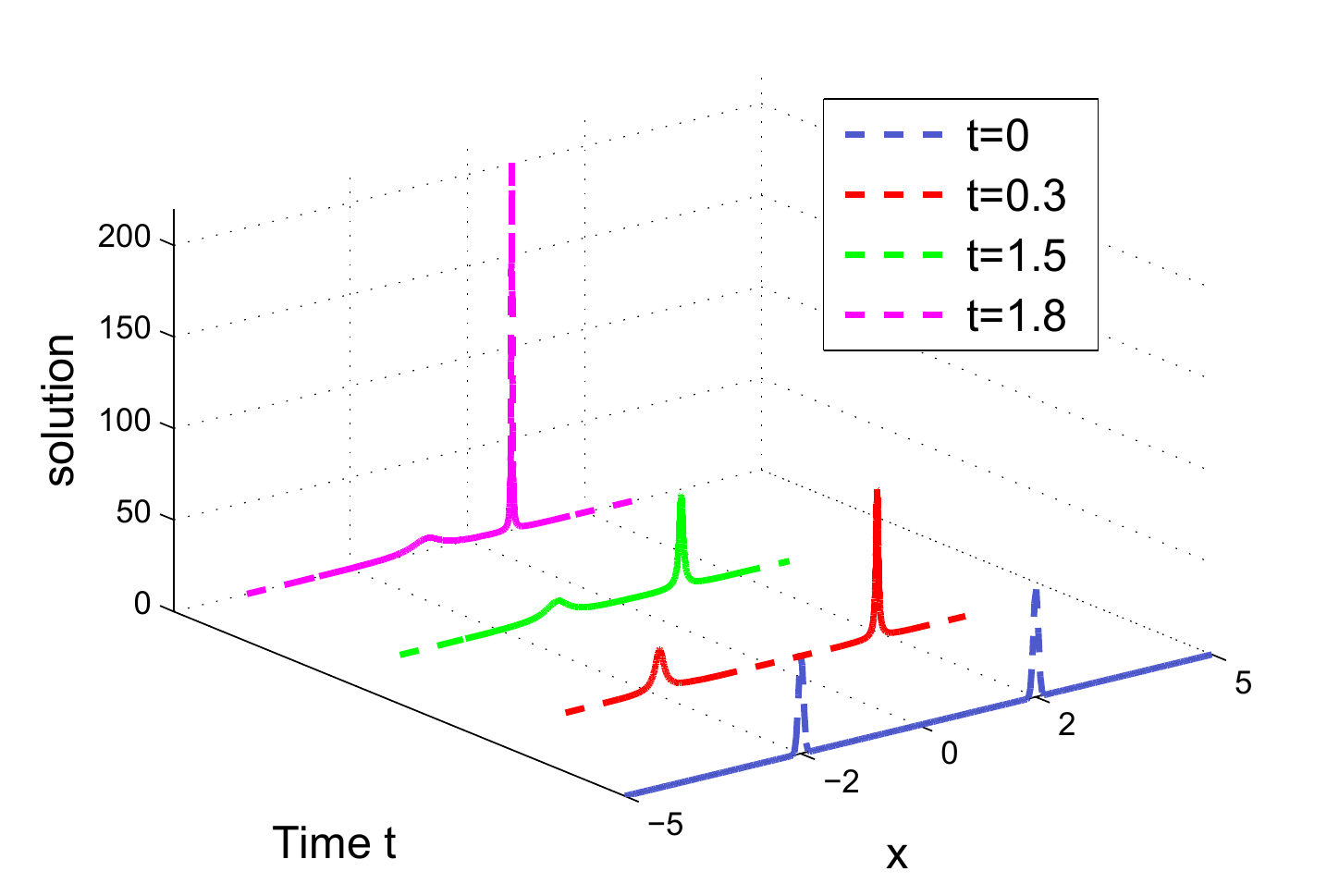}}\hspace*{0.5cm}
\subfloat{\includegraphics[width=0.42\textwidth]{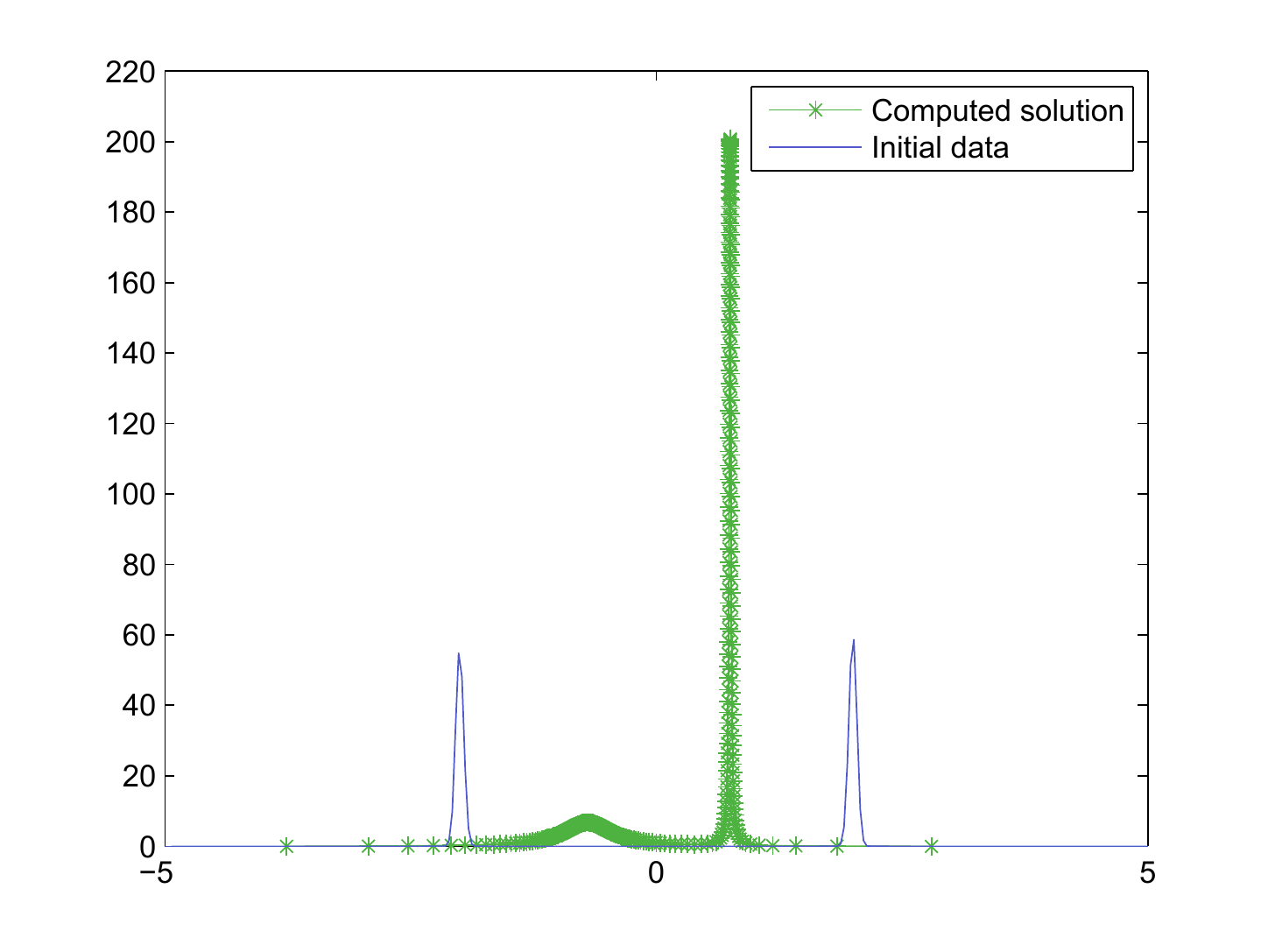}}
\caption{Evolution of the density $\rho$ for $\chi=1$ with initial datum $\rho_0 = \rho_0^1(x)$ up to time $t=1.8$.}\label{f:ex11b}
\end{center}
\end{figure}

\begin{figure}[h!]
\begin{center}
\subfloat{\includegraphics[width=0.42\textwidth]{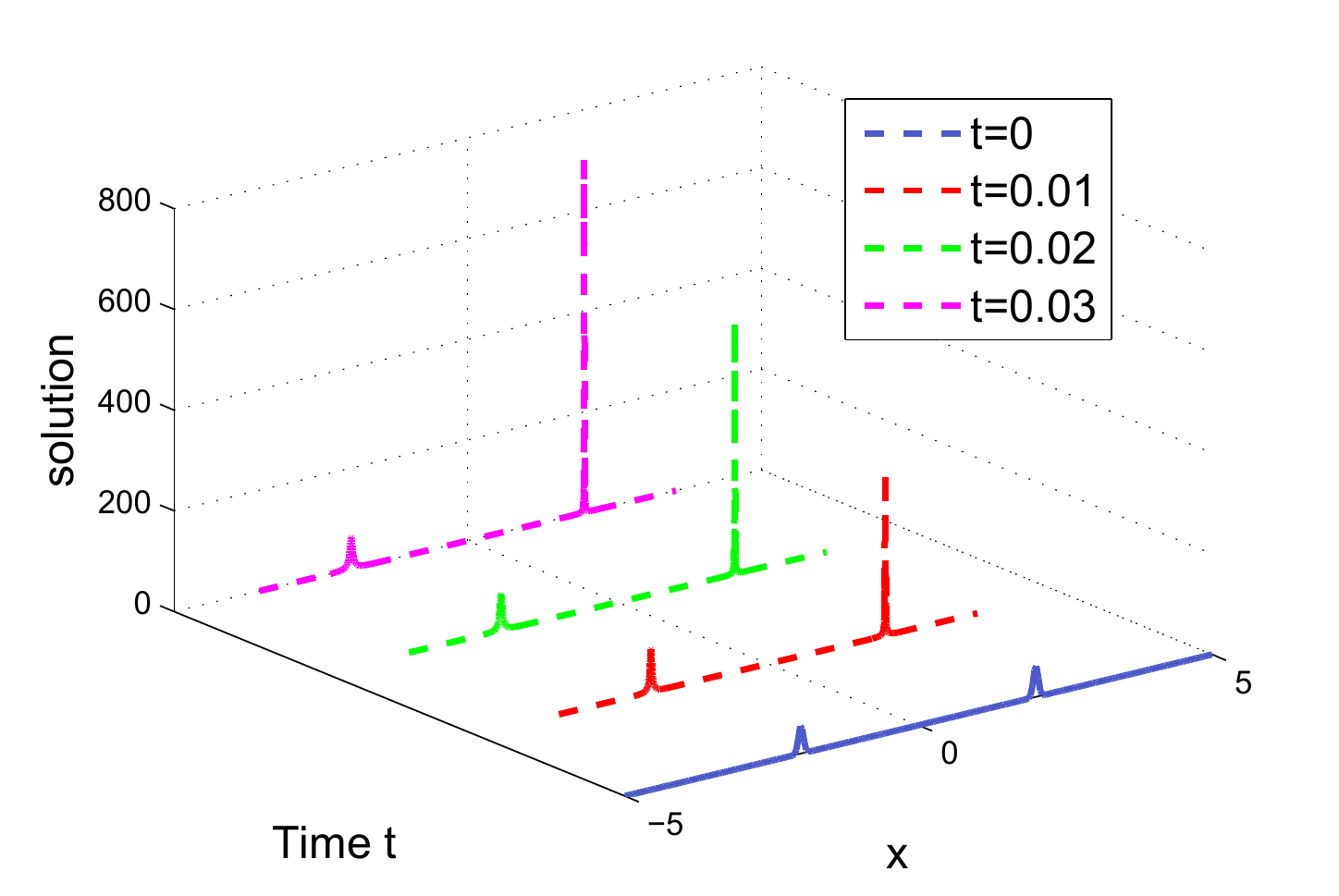}}\hspace*{0.5cm}
\subfloat{\includegraphics[width=0.42\textwidth]{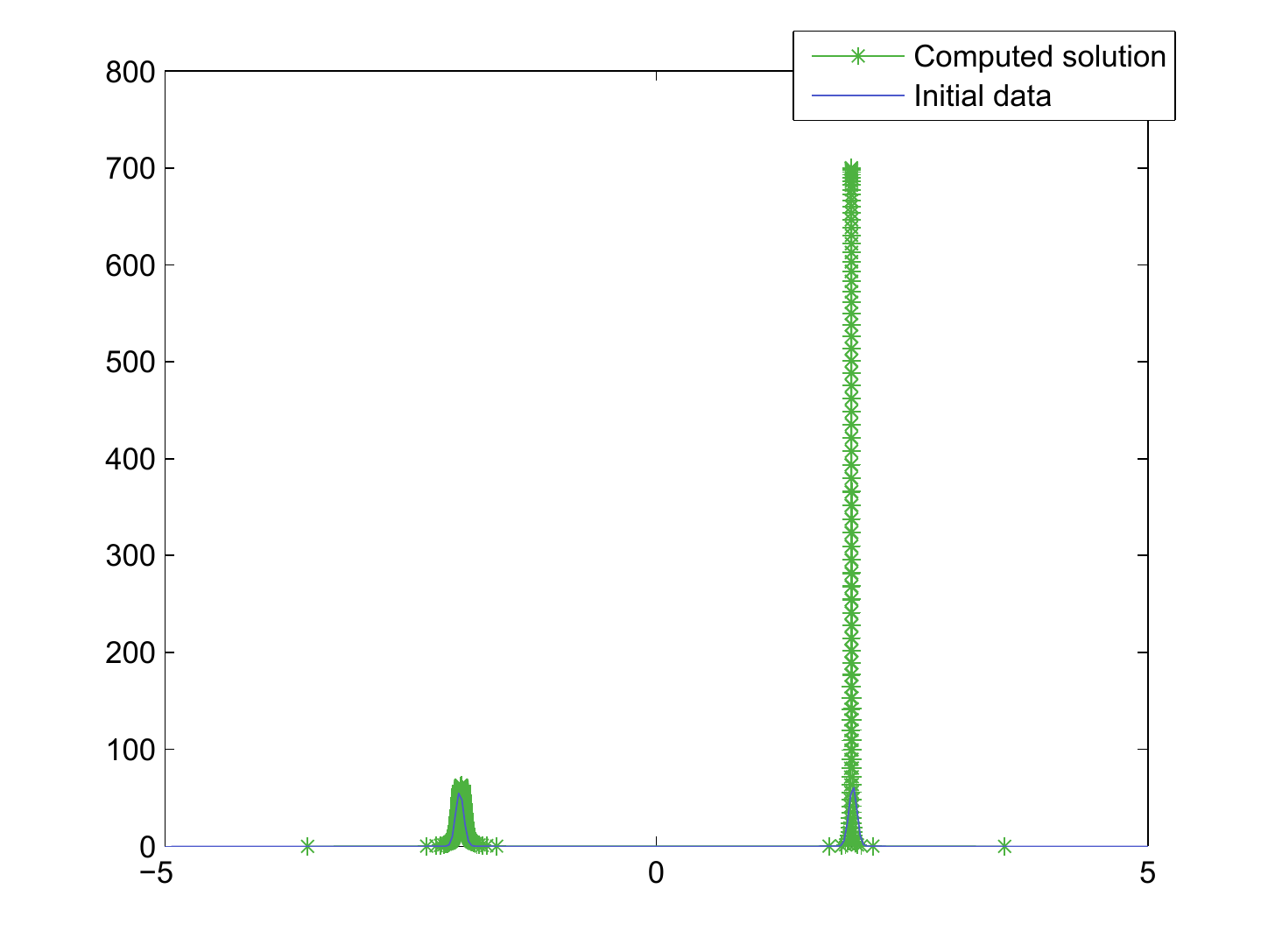}}
\caption{Evolution of the density $\rho$ for $\chi=1$ with initial datum $\rho_0 = \rho_0^2(x)$ up to time $t=0.03$.}\label{f:ex12}
\end{center}
\end{figure}

Next we study the behavior of \eqref{e:mksm} for two initial densities corresponding to the sum of two Gaussians, i.e. 
\[\rho_{0}^1(x)=\sqrt{600}\left(\frac{2\pi-0.1}{\sqrt{2\pi}}e^{\frac{-600(x-2)^2}{2}}+\frac{2\pi-0.5}{\sqrt{2\pi}}e^{\frac{-600(x+2)^2}{2}}\right),\]
and
\[\rho_{0}^2(x)=\sqrt{600}\left(\frac{2\pi+0.1}{\sqrt{2\pi}}e^{\frac{-600(x-2)^2}{2}}+\frac{2\pi-0.5}{\sqrt{2\pi}}e^{\frac{-600(x+2)^2}{2}}\right),\]
on the computational domain $\Omega = [-5,5]$ discretized into $301$ intervals of size $\Delta x = \frac{10}{301}$ with $\Delta t=10^{-3}$ and $\Delta t=10^{-4}$, respectively. The total mass is $4\pi-0.6$ in the first simulation and $4\pi-0.4$ in the second one, then they correspond to supercritical masses in which the blow-up will eventually happen. Figure \ref{f:ex11b} illustrates that if both peaks have initial masses smaller than the critical one, they initially diffuse while moving towards each other until they accumulate enough mass to blow up at the center of mass. Note that the blow-up time is not included in the figure. But if one of them has an initial mass above the critical value, the blow up happens in the center of mass of the corresponding peak before getting closer to the center of mass, see Figure \ref{f:ex12}, again the final time does not correspond to the blow-up time.

In case $M>M_c$, it is conjectured that the first blow-up should happen by the formation of a Dirac Delta with exactly the critical mass $M_c=2\pi$. Moreover, a detailed asymptotic expansion near the blow-up time for the blow-up profile was performed in \cite{HV}. The profile is not blowing-up selfsimilarly and its validity was numerically checked in \cite{BCR} for two dimensions in a numerical {\it tour-de-force}. Checking the blowup profile in time for the Keller-Segel model is quite challenging for us too since the density is a derived quantity only obtained from the diffeomorphism through \eqref{e:transf}-\eqref{e:regvarchange} by numerical differentiation. Therefore, numerical errors are likely to happen when reconstructing the density profile. 

However, we can now compare our computed diffeomorphism to the diffeomorphism representing the zeroth order term in the expansion studied by \cite{HV}. This expansion does not carry over directly to our modified KS system \eqref{e:mksm}. Nevertheless, it is easy to check following \cite{HV} that the zeroth order term should be given by a properly scaled stationary state of the critical mass case. More precisely, the solution near the blow-up time should behave like
$$
\tilde \rho (t) = \frac{2\gamma(t)^{1/2}}{1+\gamma(t) x^2} \qquad \mbox{with } \gamma(t)\simeq \rho(t,0)^2
$$
in an interval of length $2L(t)$ centered at the point of blow-up with $L(t)^{-1}\simeq \rho(t,0)$ according to \cite{HV,BCR}. In order to check this behavior, we try to avoid comparisons of the densities due to the numerical errors mentioned above and we concentrate in estimating the error between the diffeomorphisms $\Phi(t)$ and $\Psi(t)$ representing $\rho (t)$ and $\tilde\rho (t)$ respectively. The last one is computed by calculating the diffeomorphism representing the steady solution to the critical mass case $2(1+x^2)^{-1}$ and then multiplying by the dilation factor $\gamma(t)^{-1/2}$.

\begin{figure}
\subfloat[Evolution of the error in $d_W$ between $\rho(t)$ and $\tilde\rho(t)$ near the blow-up time.]{\protect\includegraphics[scale=0.5]{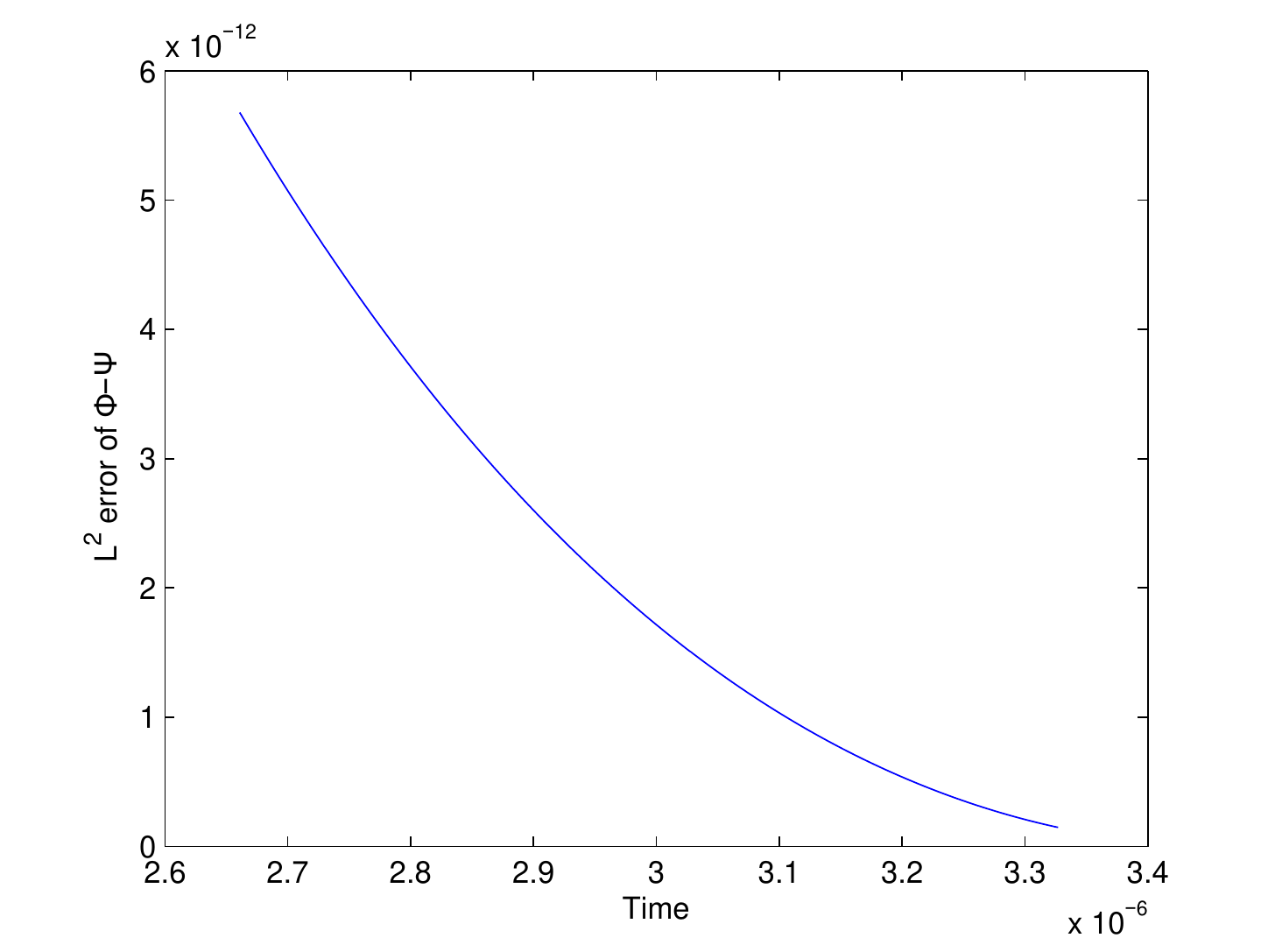}}
\subfloat[Comparison of the diffeomorphisms $\Phi$ and $\Psi$ at the blow-up time $T$. Inset shows a zoom at the center of mass.]{
\protect\includegraphics[scale=0.5]{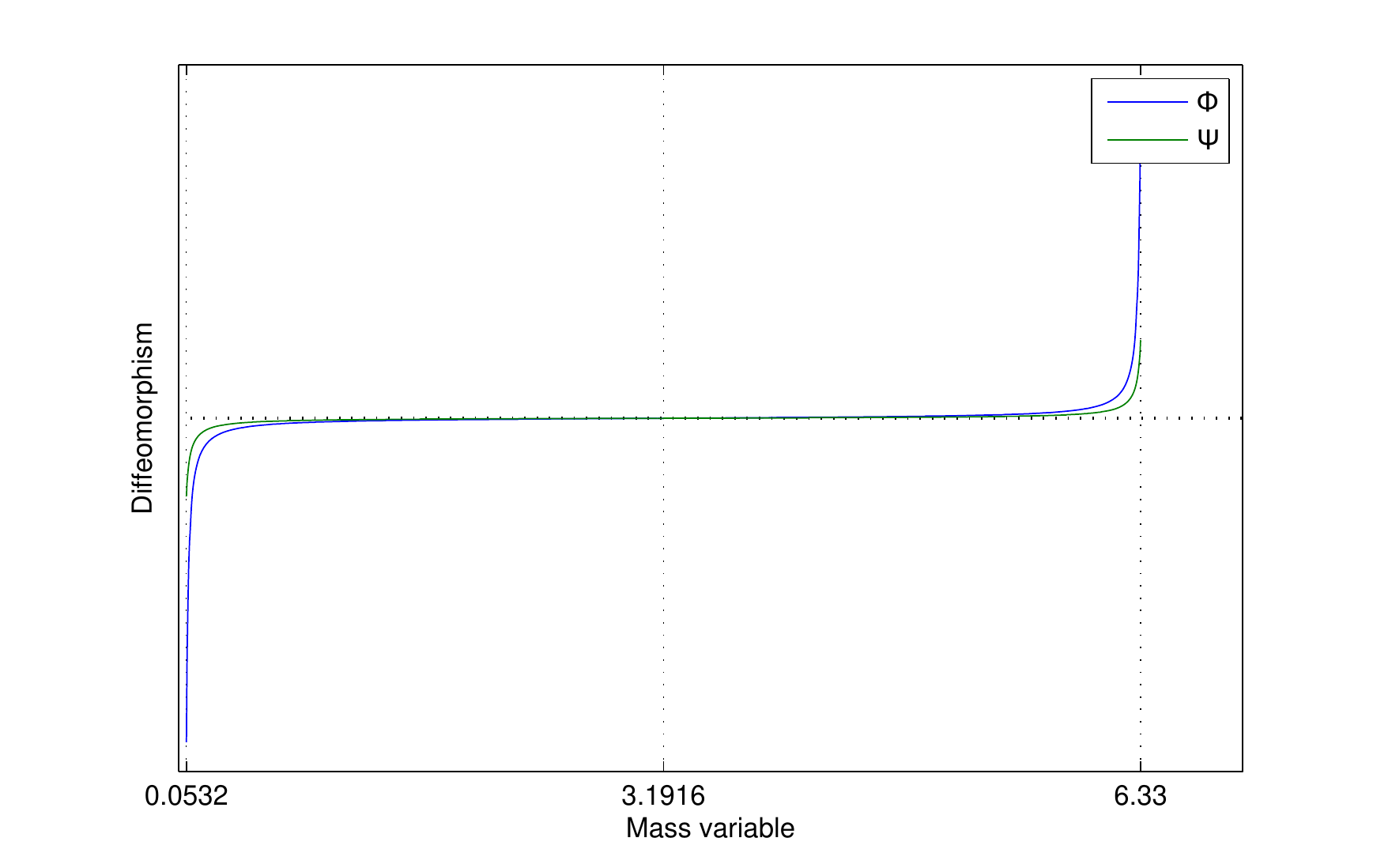}
\llap{\shortstack{%
        \includegraphics[scale=.32]{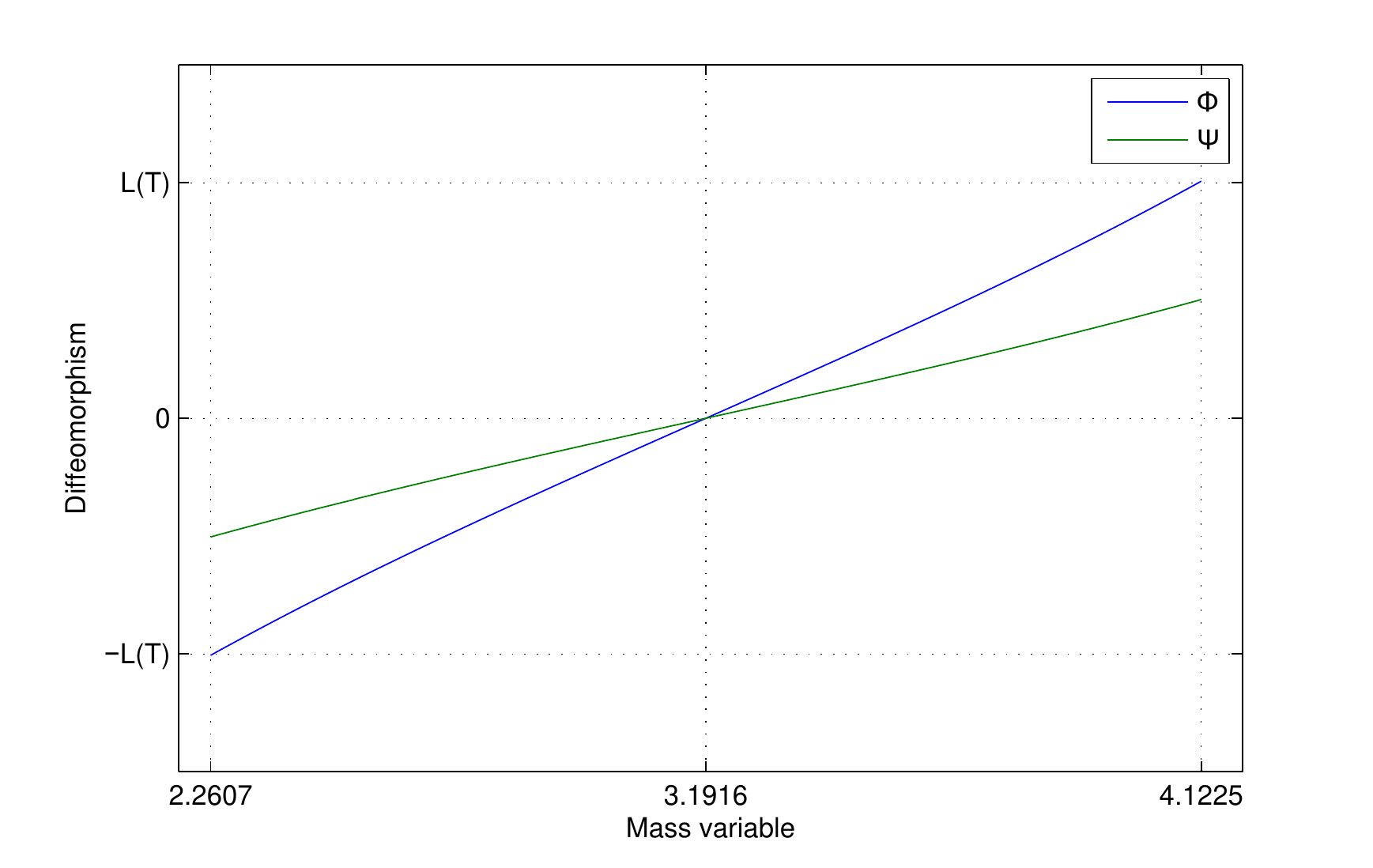}\\
        \rule{0ex}{1.25in}%
      }
  \rule{.72in}{0ex}}
}
\protect\caption{\label{fig:diffeo} Comparison of the evolution of $\Phi(t)$ and $\Psi(t)$.}
\end{figure}

We take as initial data a Gaussian already very picked at the origin, i.e. \eqref{e:gaussian} with $\sigma^2=5 \times 10^{-6}$ and $M=2\pi+0.1$ for which we have numerical blow-up quite quickly. Let $\Omega = [-0.1,0.1]$ be the computational domain discretized into $1200$ intervals of size $\Delta x = \frac{0.2}{1200}$. The discrete time steps are set to $\Delta t=10^{-9}$. Then, we compute the blow-up time $T=3.328e-6$ as the one for which our Newton's iterations do not converge. Once the blow-up time is approximated, we compute the errors in Wasserstein distance between $\rho(t)$ and $\tilde\rho (t)$ on the interval $(-L(t),L(t))$ and we can compare the corresponding diffeomorphisms at the blow-up time. Figure \eqref{fig:diffeo}(a) shows the evolution in time of the $L^2$-difference between $\Phi(t)$ and $\Psi(t)$ over the mass interval corresponding to $(-L(t),L(t))$ close to blow-up time. Figure \eqref{fig:diffeo}(b) shows the comparison between $\Phi(T)$ and $\Psi(T)$ both over the mass interval of length $2\pi$ and the zoom over the mass interval corresponding to $(-L(T),L(T))$.

\subsection{2D simulations}

\noindent Next we present several numerical simulations in 2D. Let $\Omega$ denote the unit circle, $\mathcal{M} = \bigcup_{V \in \mathcal{V}} x_V$ 
the computational mesh with vertexes $x_V$, $\mathcal{V} = \lbrace 1, \ldots n_V\rbrace$ and $\mathcal{N} = \bigcup_{V \in \mathcal{V}}\Phi(x_V)$ the transformed mesh. We choose the following parameters
if not stated otherwise:
\begin{itemize}
\item $\Delta t=0.01$
\item $n_V=1758$ which corresponds to $3386$ triangles of maximum size $h_{\max} = 0.05$.
\end{itemize}
The initial datum is given by Gaussian centered around the origin
\begin{align*}
\rho_0(x) = \frac{1}{2\pi \sigma^2} e^{-\frac{1}{2\sigma^2}(x^2+y^2)}+c \text{ with } \sigma = 0.3,
\end{align*}
where $c$ is chosen such that $\int_{\Omega} \rho_0(x) dx = 1$. In the pre-processing step we solve the heat equation on the time interval $t \in (0,2]$ using an implicit in time discretization and $H^1$-conforming elements of order 6. The solution is computed at $t_i=0.002i,\,i=1,\ldots,1000$ and used for the computation of the initial diffeomorphism. The error bounds for the Newton schemes are set to $\epsilon_1 = \epsilon_2 = 10^{-6}$, the
regularization parameter $\varepsilon$ in the post-processing step to $\varepsilon = 10^{-2}$. All solutions are computed as detailed in Algorithm \ref{alg}.

\subsubsection{Attraction potentials}
\begin{figure}[h!]
\begin{center}
\subfloat[Transformed mesh]{\includegraphics[width=0.3\textwidth]{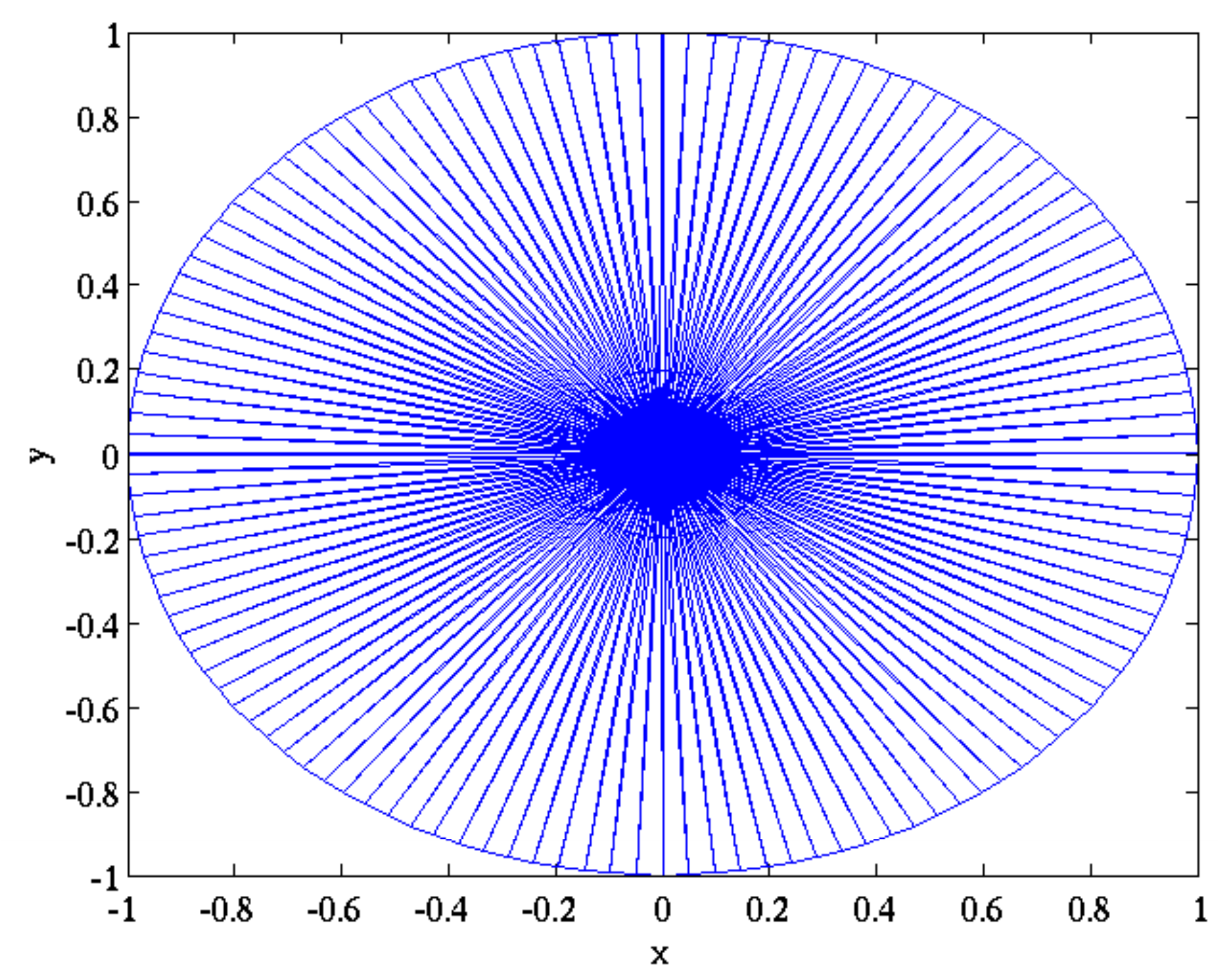}}\hspace*{0.25cm}
\subfloat[Density $\rho$]{\includegraphics[width=0.3\textwidth]{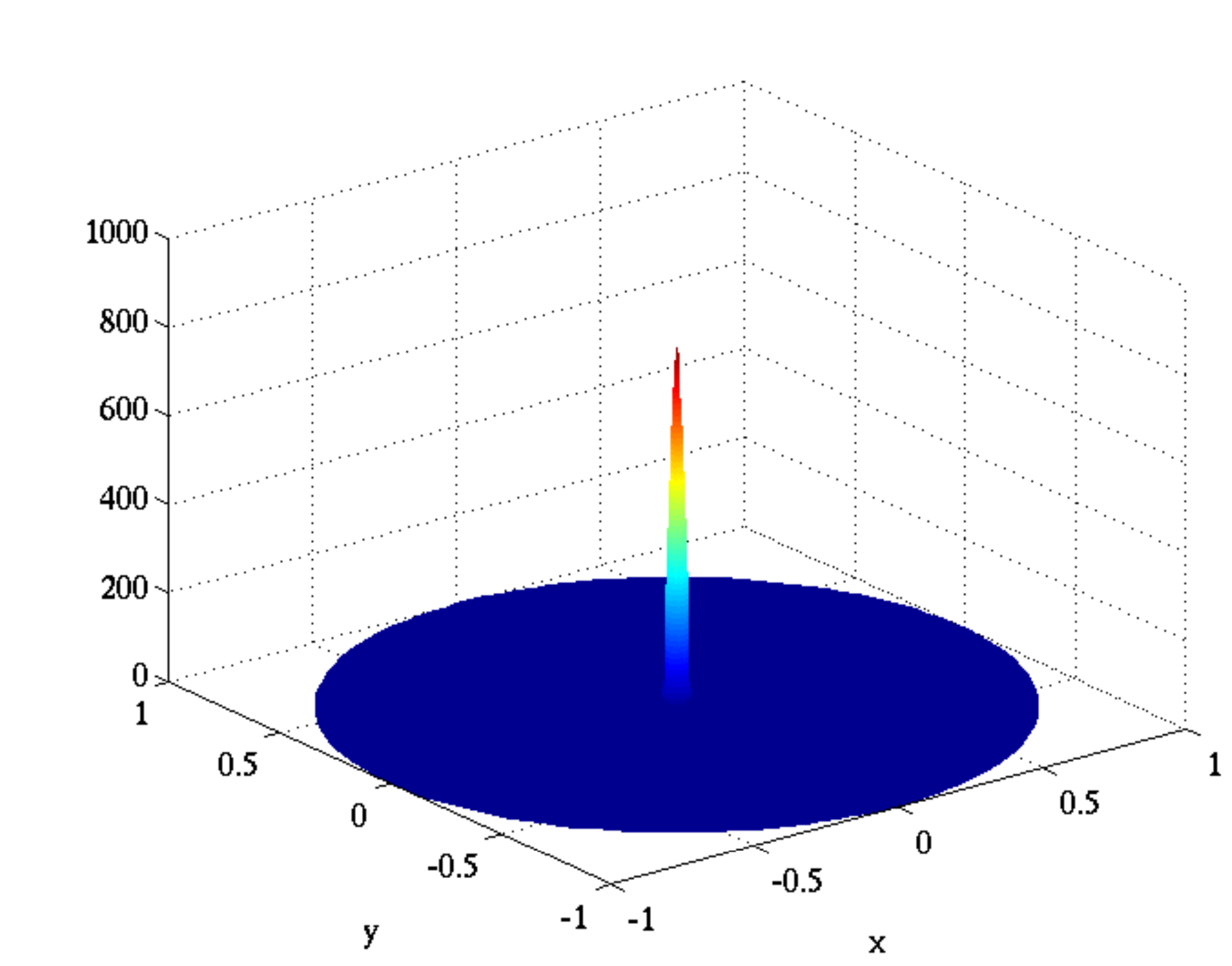}} \hspace*{0.25cm}
\subfloat[$d_W(\rho(t), \delta_0)$.]{\includegraphics[width=0.3\textwidth]{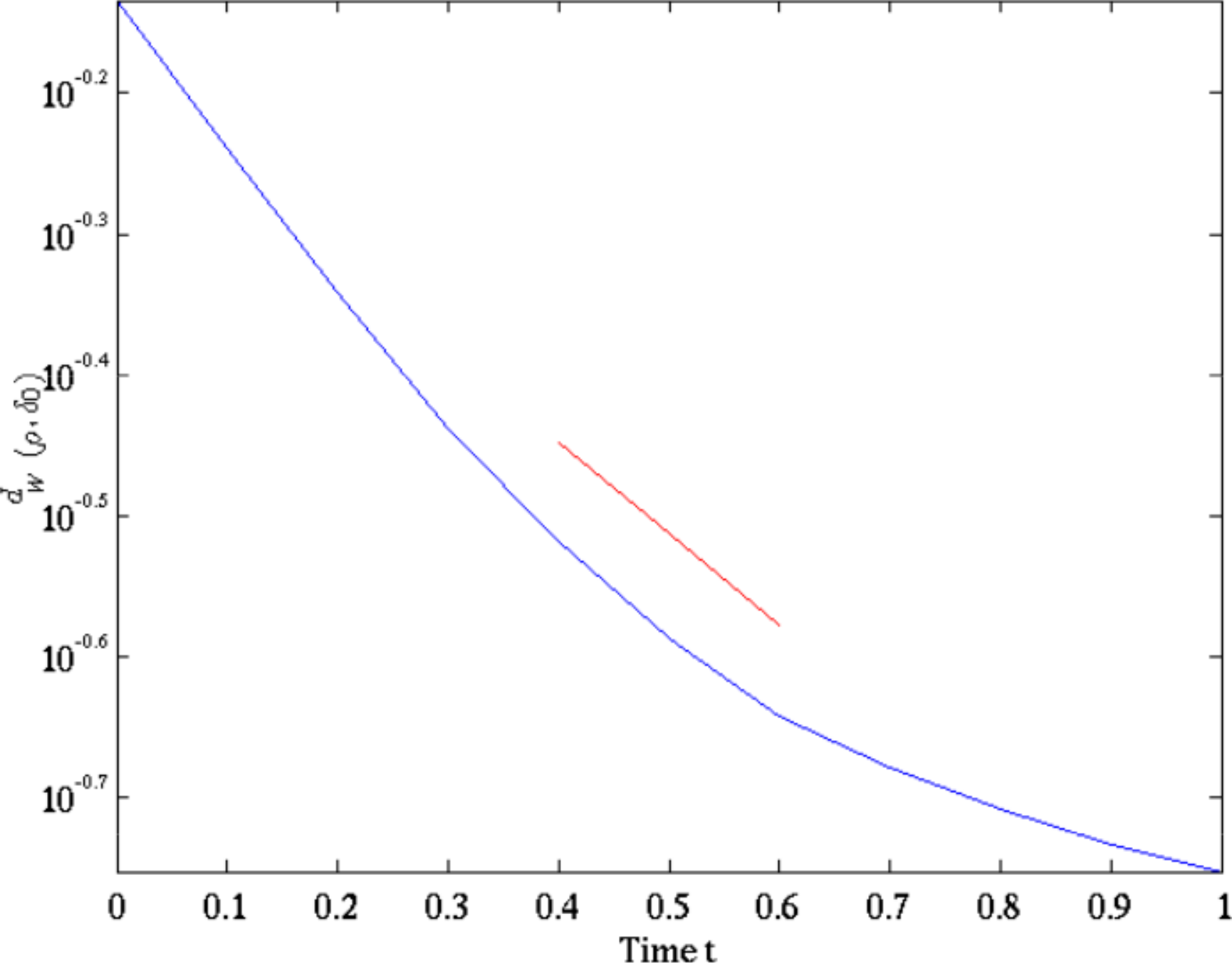}}
\caption{Simulation results in the case of a purely attractive potential $W = \frac{1}{2} \lvert x \rvert^2$.} \label{f:ex1}
\end{center}
\end{figure}

\noindent First we consider an attraction potential of harmonic type
$W(x) = \frac{1}{2} \lvert x \rvert^2$.
In the case of a purely attractive potential we expect the formation of a Delta Dirac at the center of the domain. Figure \ref{f:ex1} illustrates the formation of the blow up. In the first picture we see the transformed mesh, where all triangles are moved towards the
center of the domain. The second and third picture show the reconstructed density profile $\rho$ as well as the decay of the distance towards the Dirac Delta in time. The numerical simulations confirm the theoretical results, indicated by the red line of slope $-1$, at the beginning of the simulation. The bad match towards the end  results most likely from the formation of the Delta Dirac and the inaccuracy caused by it.

\subsubsection{Attraction-repulsion potentials}

\noindent Next consider the attraction-repulsion potential \eqref{e:W2}. In the case $a=4$ and $b=2$ the solution concentrates on a ring of radius $r = \frac{1}{2}$, see \cite{BH1,BH2}.
\begin{figure}[h!]
\begin{center}
\subfloat[Transformed mesh]{\includegraphics[width=0.32\textwidth]{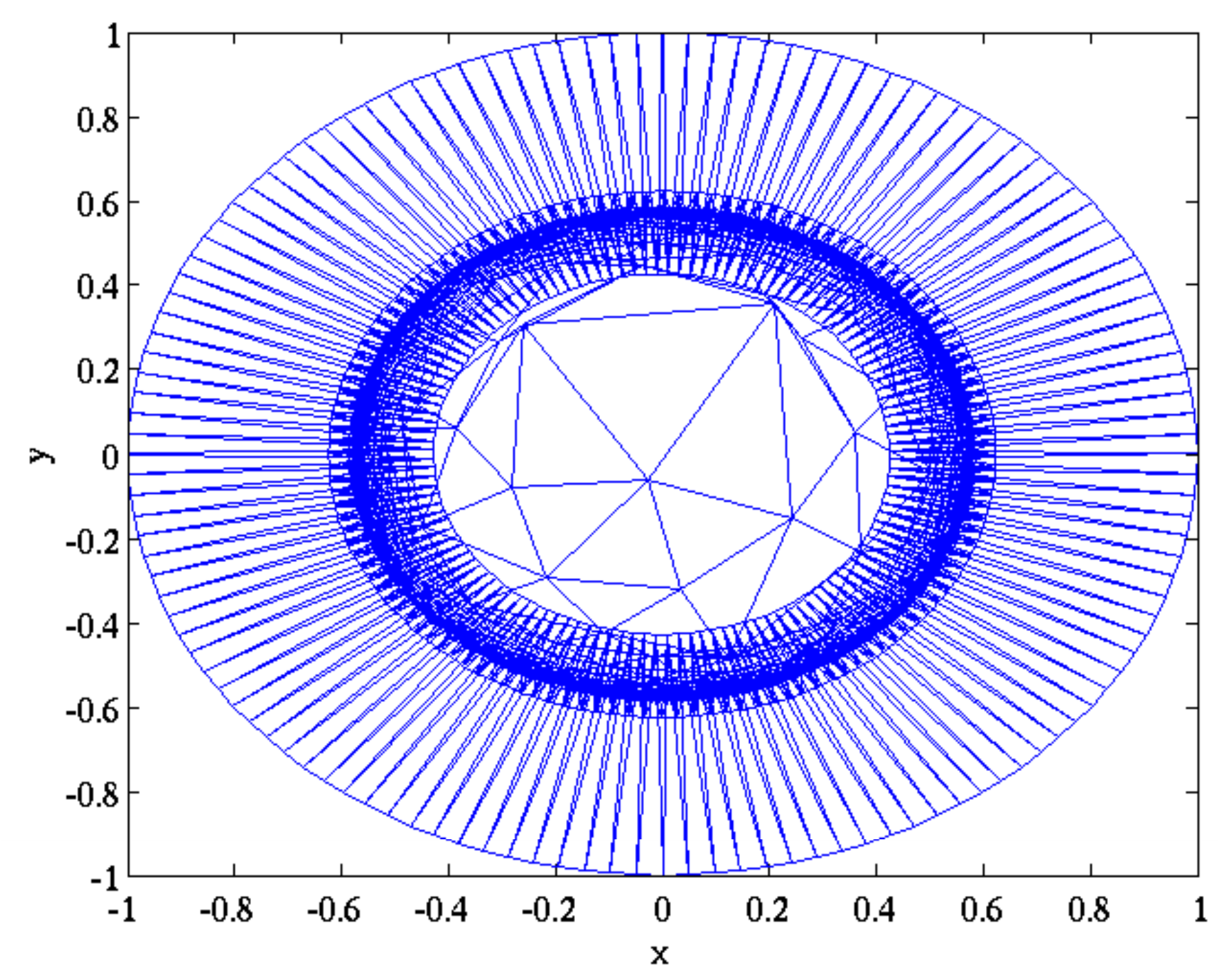}}\hspace*{0.25cm}
\subfloat[Density $\rho$]{\includegraphics[width=0.32\textwidth]{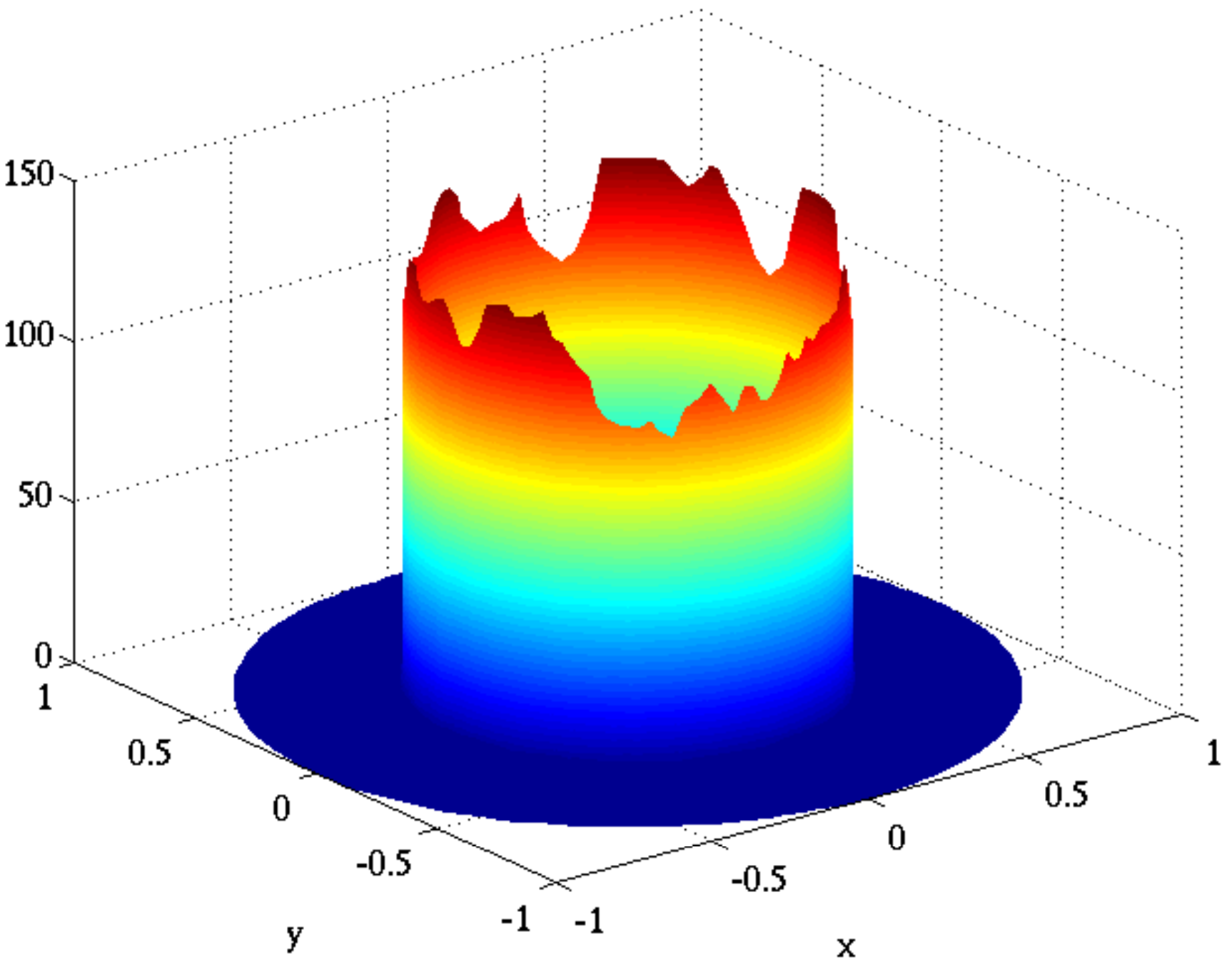}}\hspace*{0.25cm}
\subfloat[Relative entropy]{\includegraphics[width=0.32\textwidth]{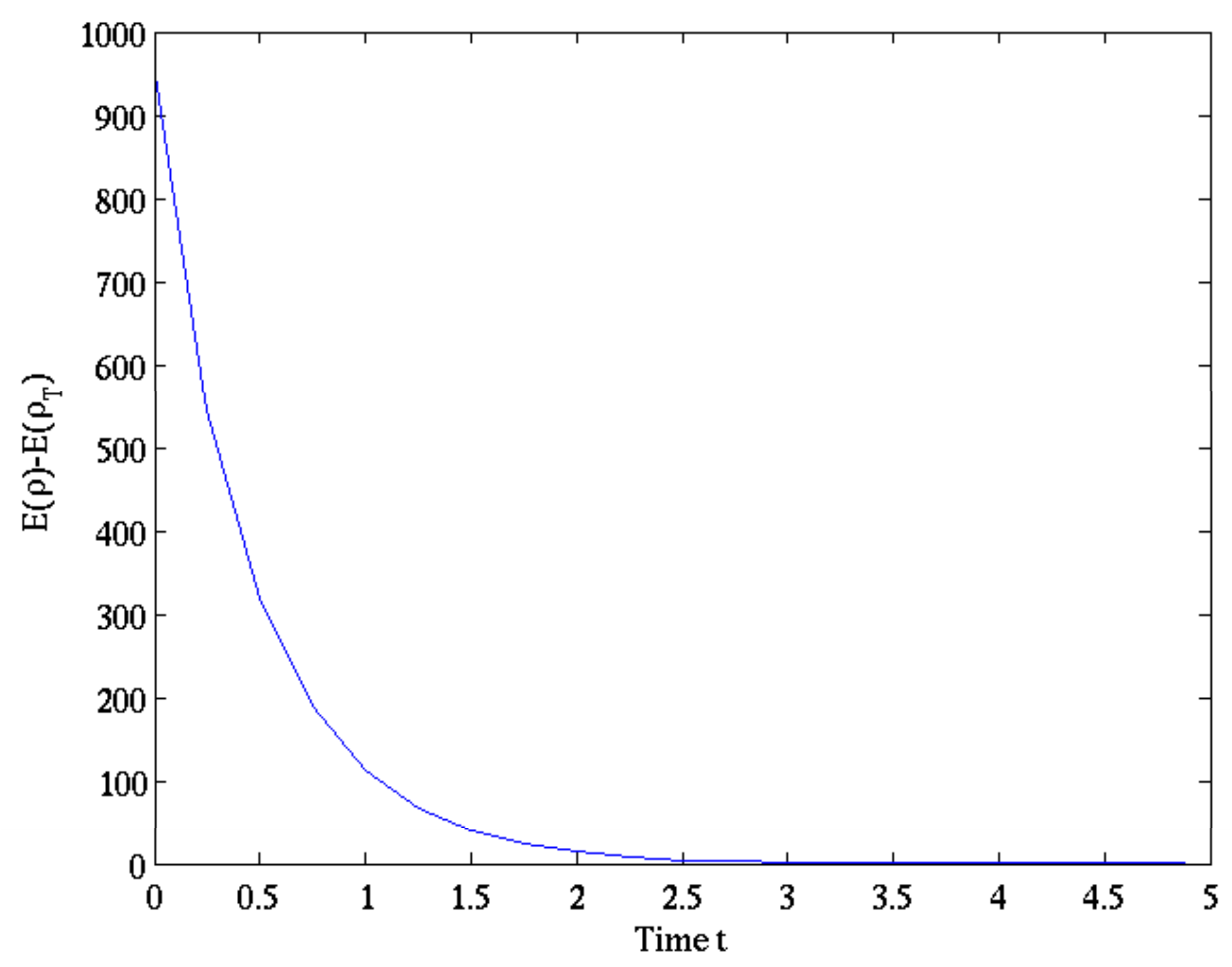}}
\caption{Simulation results for an attractive-repulsive potential  $W = \frac{1}{4} \lvert x \rvert^4 - \frac{1}{2} \lvert x \rvert^2$.} \label{f:ex2}
\end{center}
\end{figure}
The formation of the Delta Dirac ring is clearly visible in Figure \ref{f:ex2}. The first picture corresponds to the transformed initial mesh, the second and third show the reconstructed density and the decay of the relative entropy respectively. Note that we map the boundary of the computational domain onto itself, hence the boundary nodes do not move.

\begin{figure}[h!]
\begin{center}
\subfloat[Transformed mesh]{\includegraphics[width=0.32\textwidth]{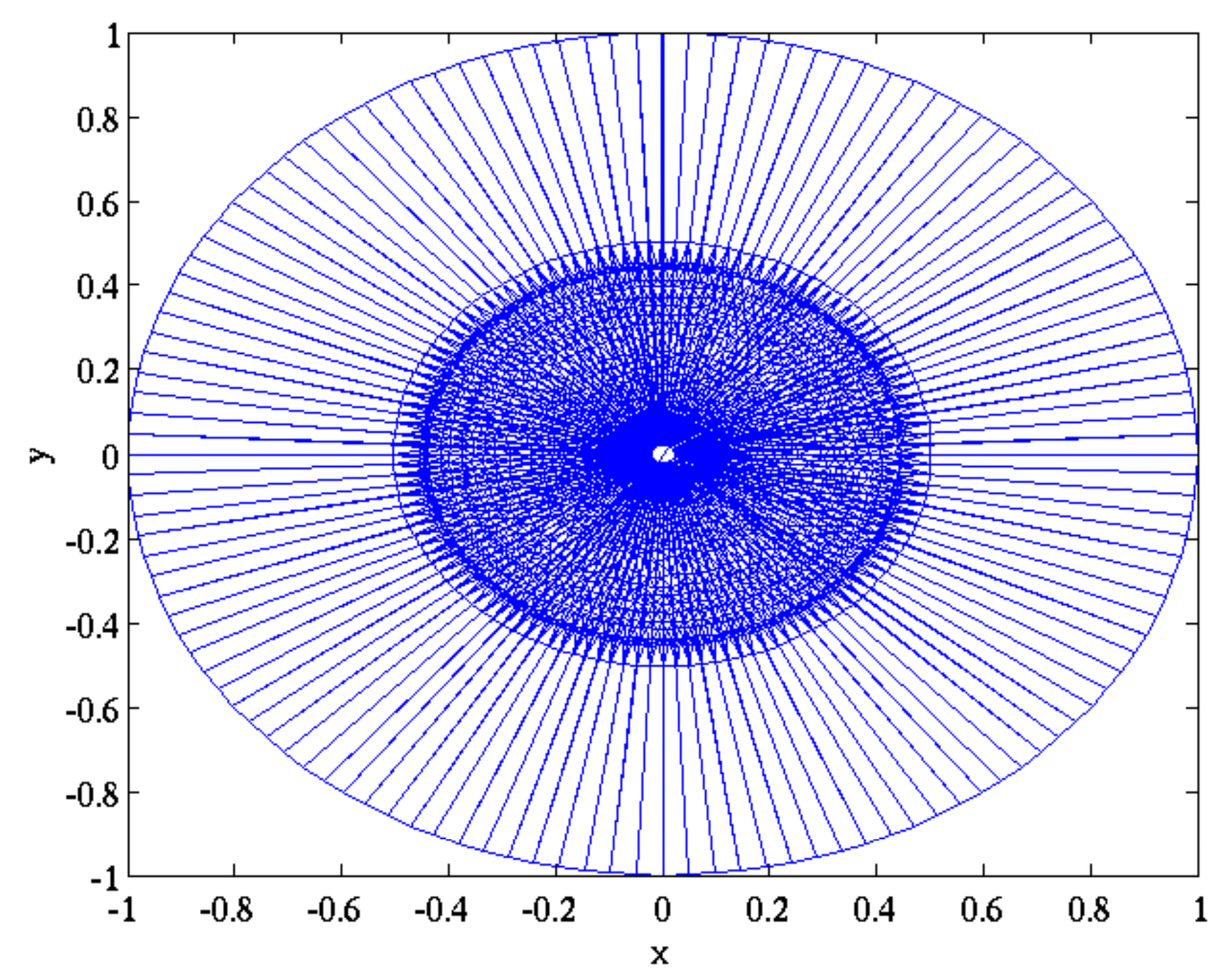}}\hspace*{0.25cm}
\subfloat[Density $\rho$]{\includegraphics[width=0.32\textwidth]{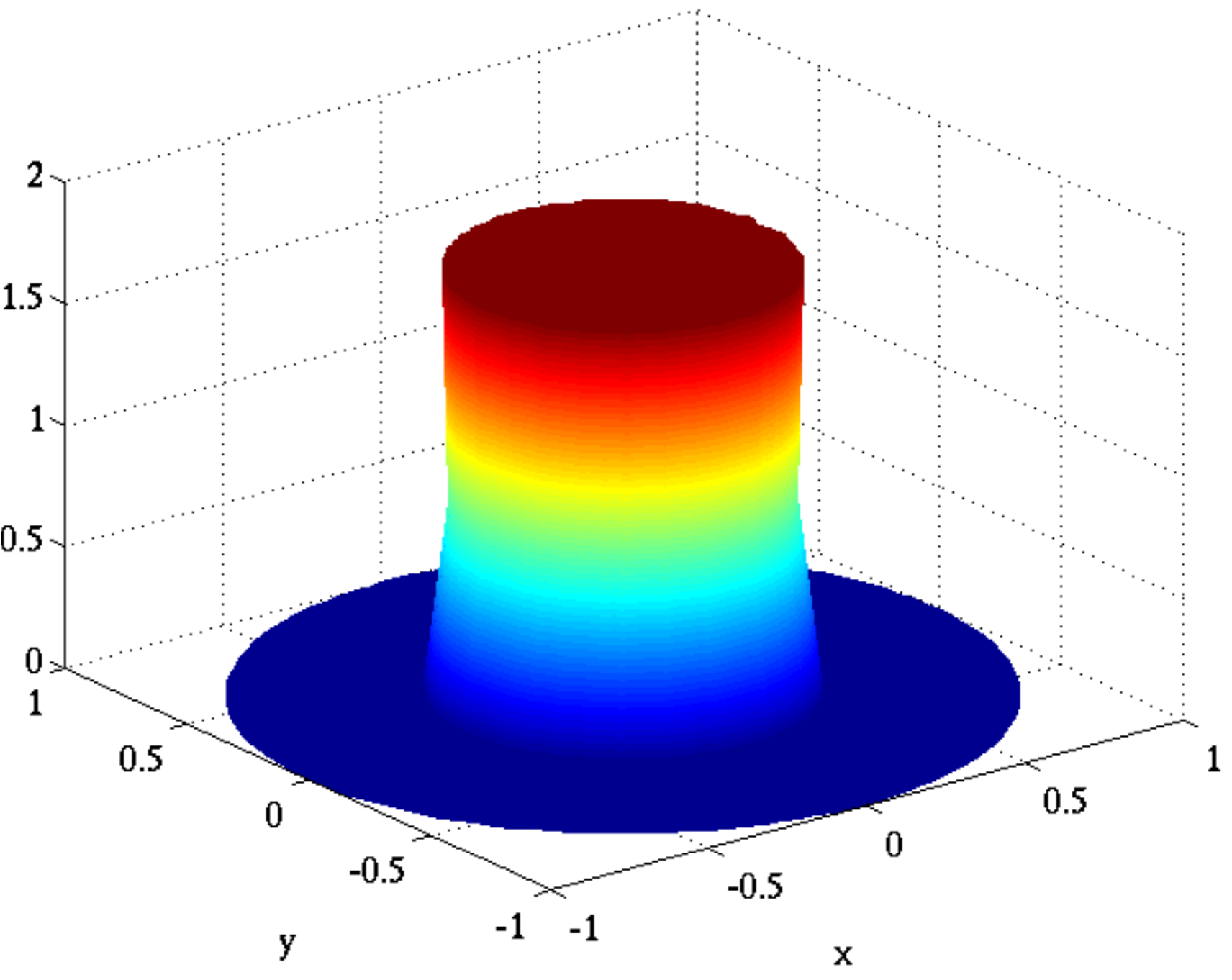}}\hspace*{0.25cm}
\subfloat[Entropy]{\includegraphics[width=0.32\textwidth]{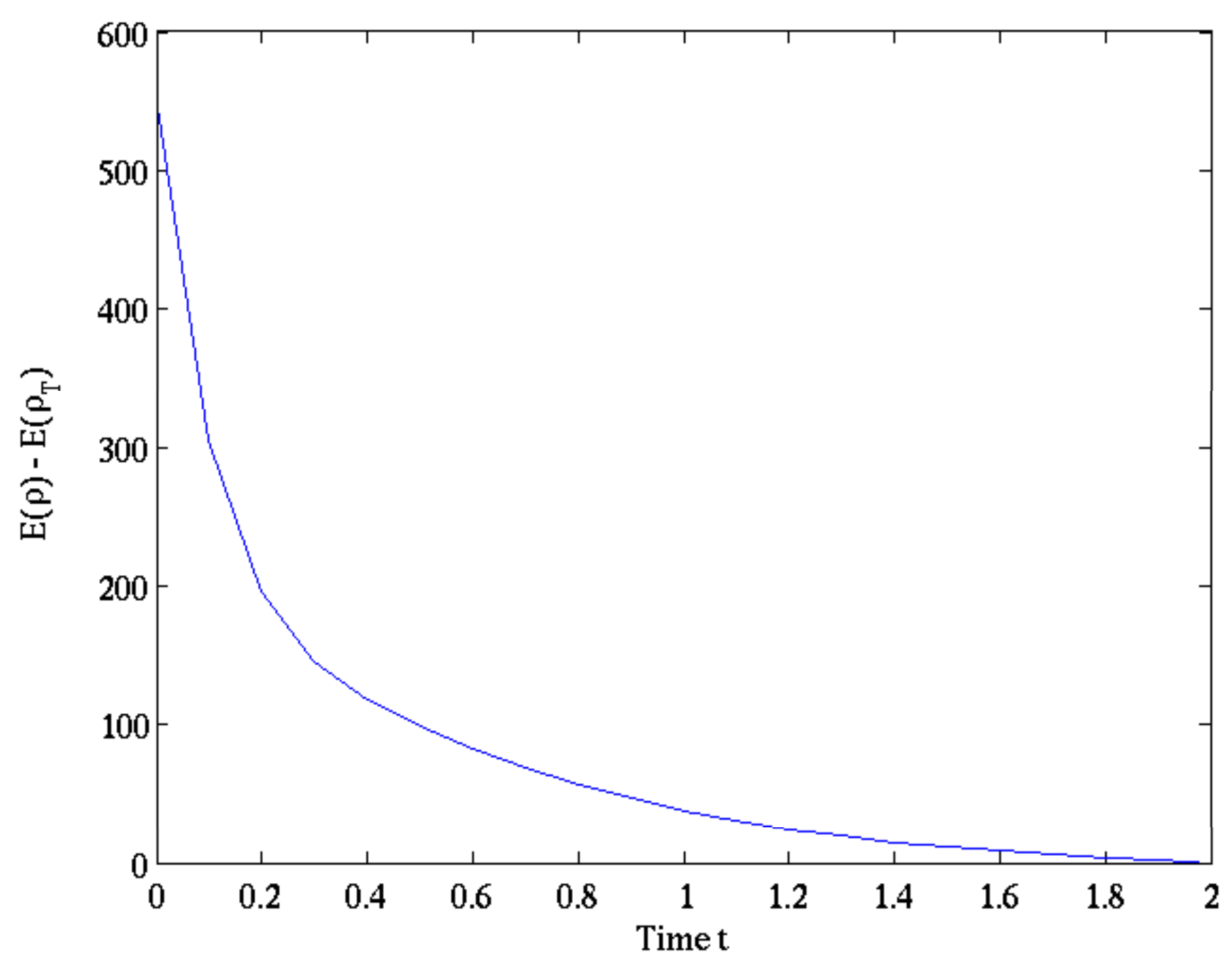}}
\caption{Simulation result for an attractive-repulsive potential $W = \frac{1}{2} \lvert x \rvert^2 - \ln(\lvert x \rvert)$.} \label{f:ex3}
\end{center}
\end{figure}

In the case of a logarithmic repulsion, that is $a = 2$ and $b=0$ in \eqref{e:W2}
we expect the formation of a compactly supported steady state circle with radius $r = \frac{1}{2}$.  We see that the vertexes of the mesh concentrate in a circle of radius $r=\frac{1}{2}$, except for the boundary nodes which are fixed, see Figure \ref{f:ex3}. Again the second and third plot correspond to the reconstructed density $\rho$ at time $t=2$ and the decay of the relative entropy functional.

\begin{figure}[h!]
\begin{center}
\subfloat[Transformed mesh]{\includegraphics[width=0.32\textwidth]{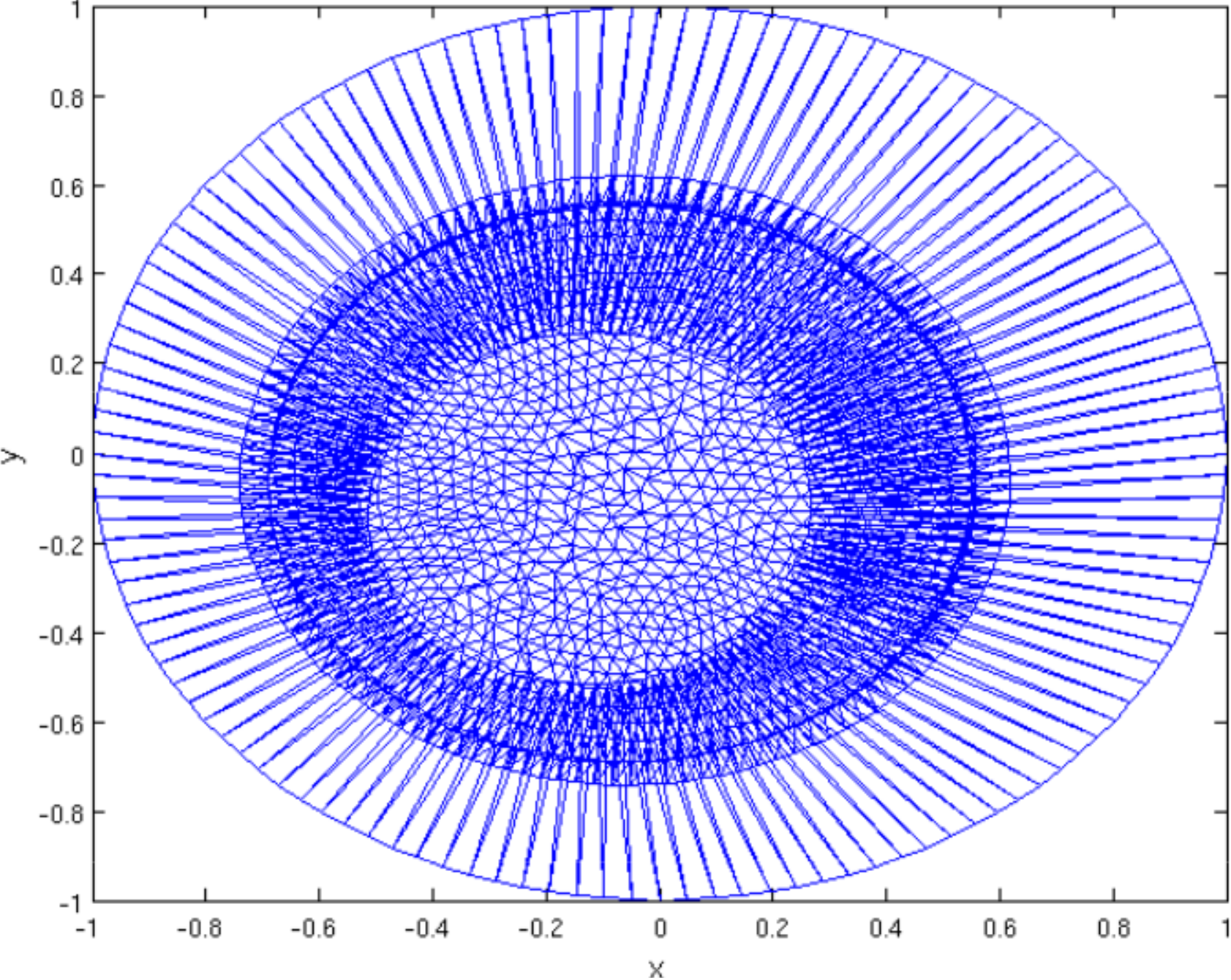}}\hspace*{0.25cm}
\subfloat[Density $\rho$]{\includegraphics[width=0.32\textwidth]{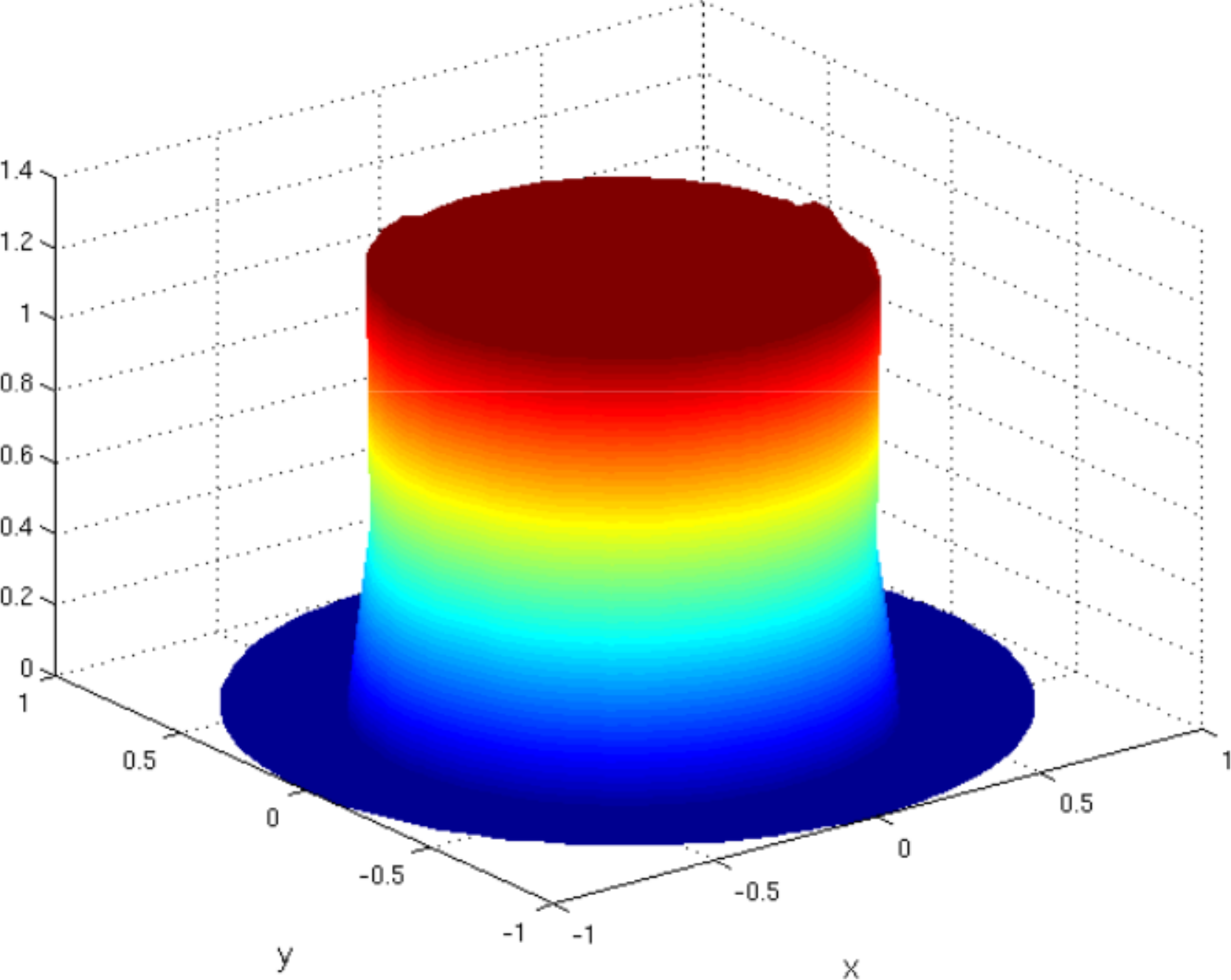}}\hspace*{0.25cm}
\subfloat[Relative entropy]{\includegraphics[width=0.32\textwidth]{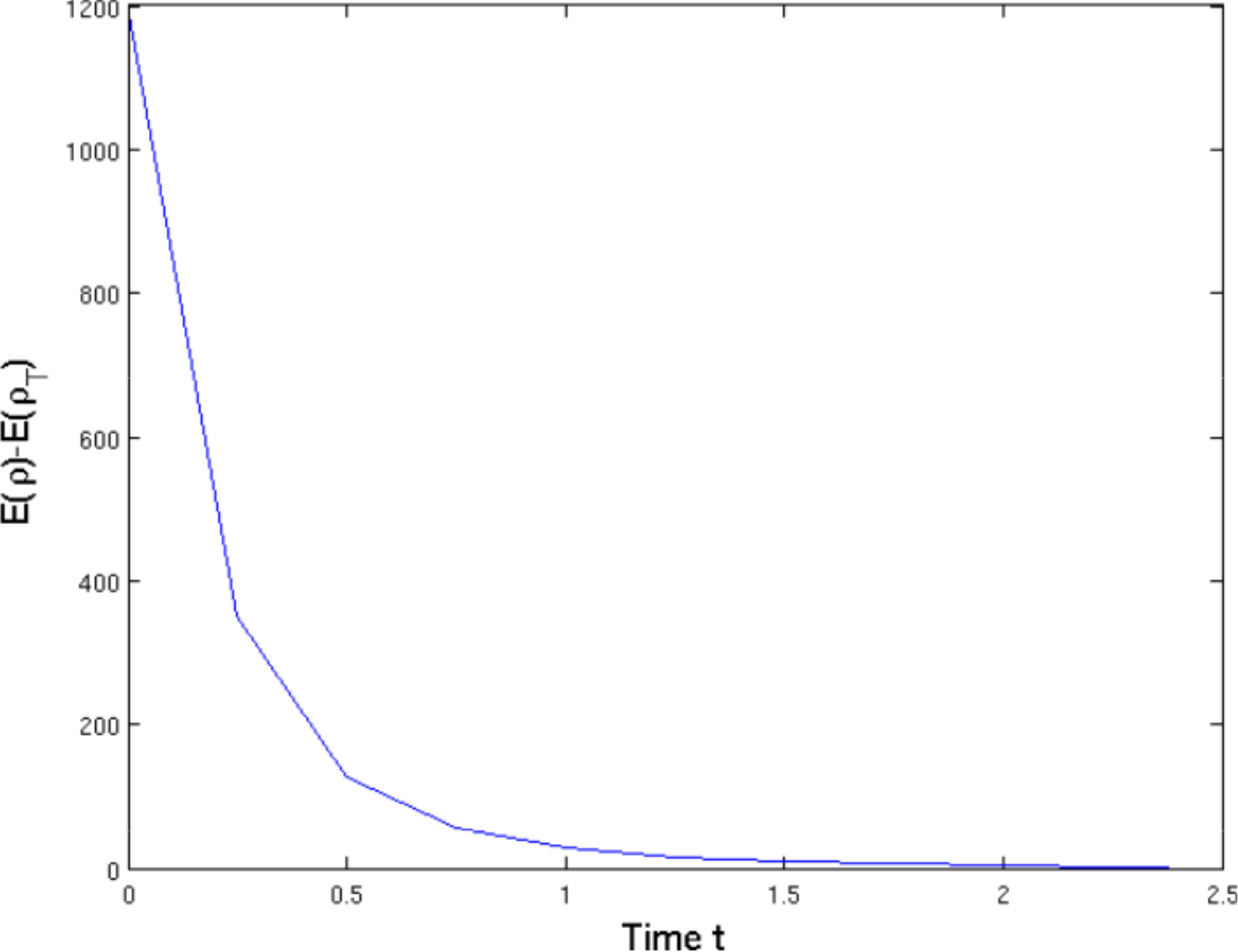}}
\caption{Simulation results for an attractive-repulsive potential $W = \frac{1}{2} \lvert x \rvert^2 - \ln(\lvert x \rvert)$ in case of a not radially symmetric initial datum.} \label{f:ex6}
\end{center}
\end{figure}

Next we consider the agreggation equation with the logarithmic repulsive potential and harmonic confinement,  $a = 2$ and $b=0$ in \eqref{e:W2}, with an additional external potential of the form 
$V(x) = -\frac{\alpha}{\beta} \ln(\lvert x \rvert)$.
\begin{figure}[h!]
\begin{center}
\subfloat[Transformed mesh]{\includegraphics[width=0.32\textwidth]{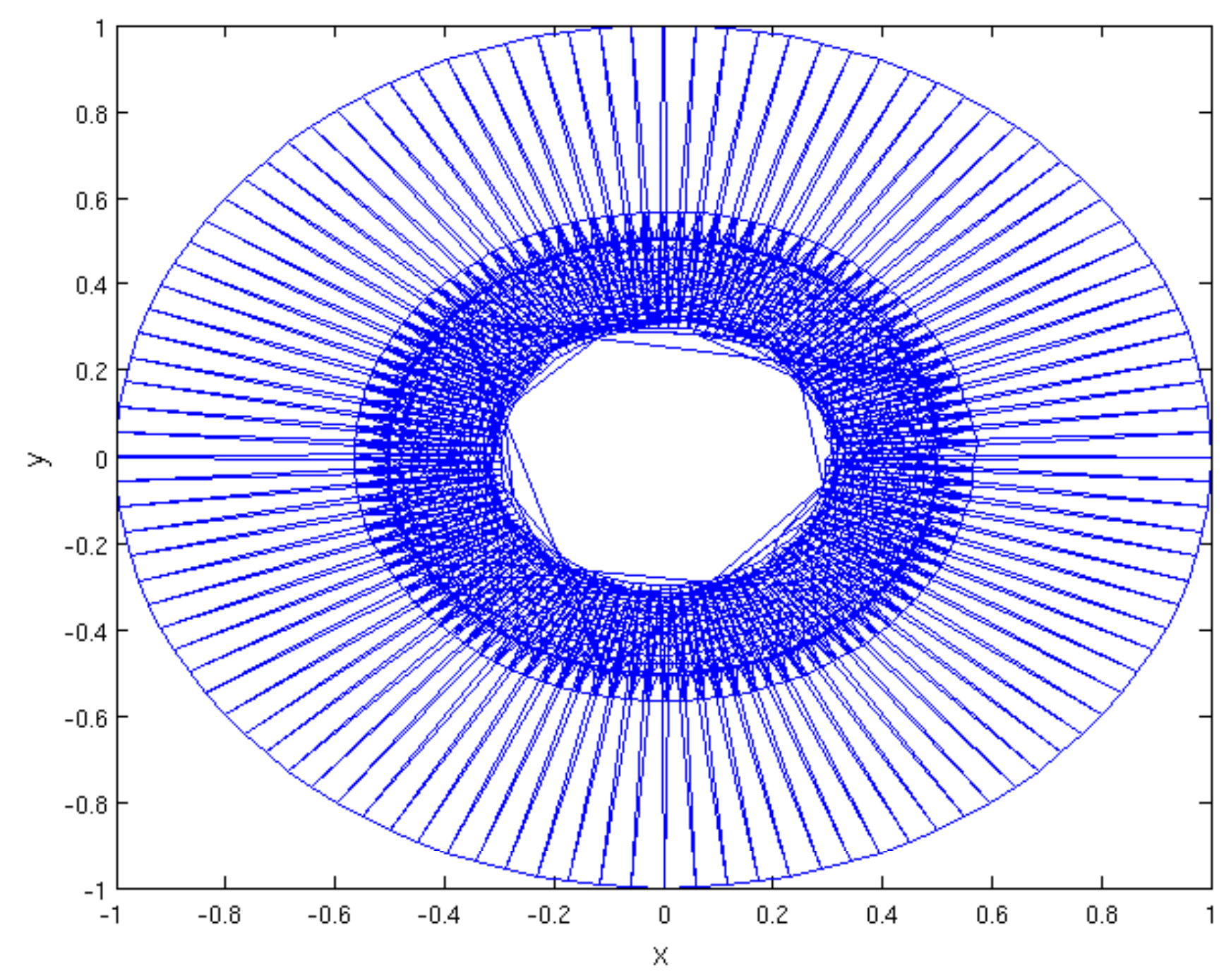}}
\subfloat[Density $\rho$]{\includegraphics[width=0.32\textwidth]{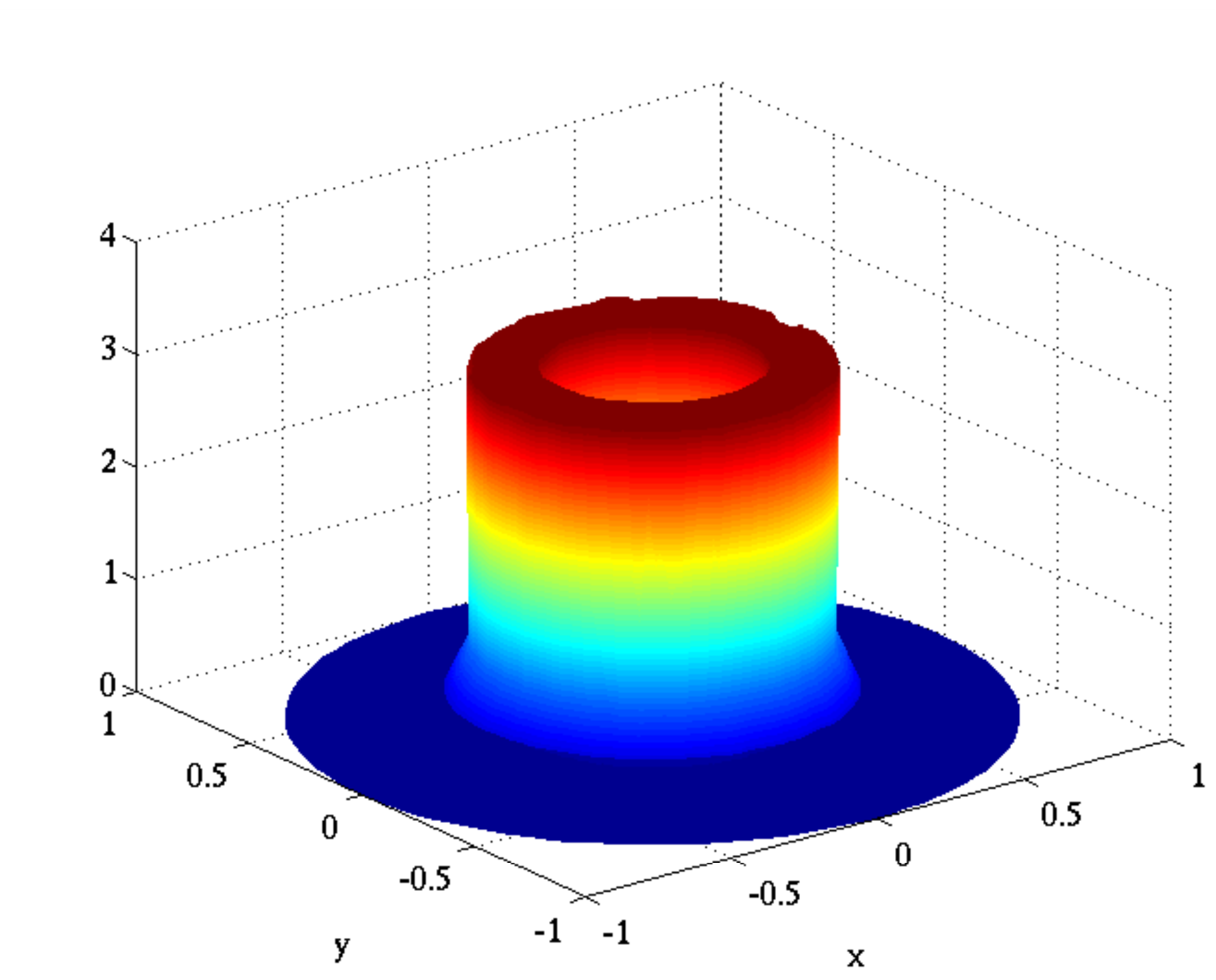}}
\subfloat[Relative entropy]{\includegraphics[width=0.32\textwidth]{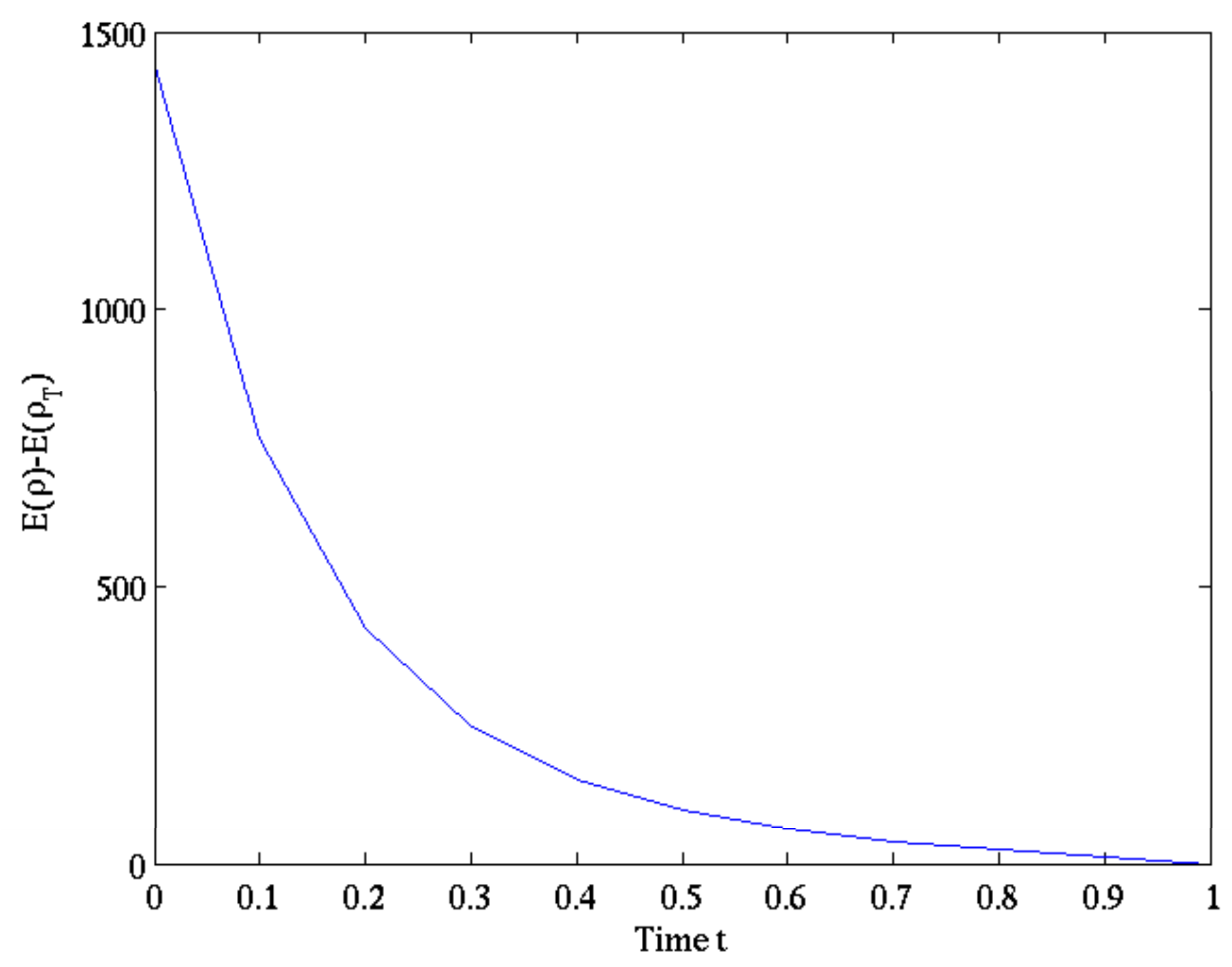}}
\caption{Simulation results for an attractive-repulsive potential $W = \frac{1}{2} \lvert x \rvert^2 - \ln(\lvert x \rvert)$ and an additional potential $V = \frac{1}{4} \ln (\lvert x \rvert)$.} \label{f:ex4}
\end{center}
\end{figure}
This model has been proposed as a way to find the spatial shape of the milling profiles in microscopic models for the dynamics of bird flocks \cite{CDP}. We obtain this steady state also in the numerical simulations, see Figure \ref{f:ex4}. The inner and outer radius of this mill depend on the parameters $\alpha$ and $\beta$ and are given
\begin{align*}
R_i = \sqrt{\frac{\alpha}{\beta}} \textrm{ and } R_o = \sqrt{\frac{\alpha}{\beta} + \frac{1}{2\pi}}.
\end{align*}
see \cite{CCH2015}. The transformed mesh clearly shows the formation of the steady state annulus with $\alpha=1$ and $\beta=4$. Note that the 'vacuum formation' at the center distorts the mesh at the center. Note that the distortion of the triangles in Figure \ref{f:ex4} does not affect the performance of the numerical method itself, since all computations are done using the original mesh on the reference domain $\tilde{\Omega}$. It only affects the reconstruction of the final density $\rho$, where we observe numerical instabilities close to the support of the steady state annulus. Therefore we choose to solve a regularized version of \eqref{e:varchange}.\\ 
If we start with a not radially symmetric initial datum of the form 
\begin{align*}
\rho_0(x) = \frac{1}{4\pi 0.25^2} e^{-\frac{1}{2 0.25^2}((x-0.15)^2(y-0.25)^2)} + \frac{1}{4\pi 0.2^2} e^{-\frac{1}{20.2^2}((x+0.3)^2+(y+0.4)^2)}
\end{align*}
and again consider a potential of the form \eqref{e:W2} with $a=2,~b=0$, we observe the formation of an off-centered compactly supported bump, see Figure \ref{f:ex6}.

\subsubsection{The Keller-Segel model}

\begin{figure}[h!]
\begin{center}
\subfloat[Transformed mesh]{\includegraphics[width=0.32\textwidth]{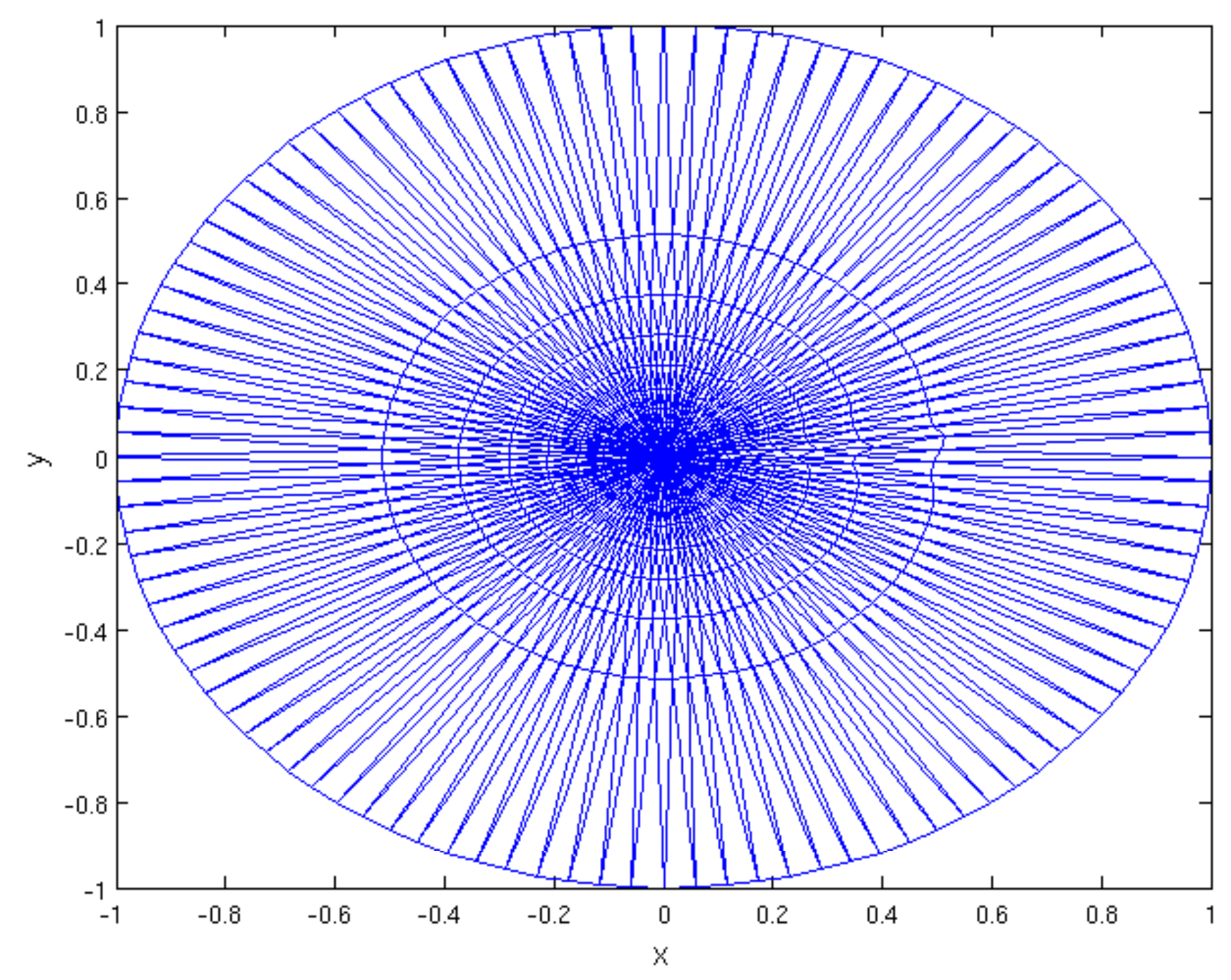}}
\subfloat[Density $\rho$]{\includegraphics[width=0.32\textwidth]{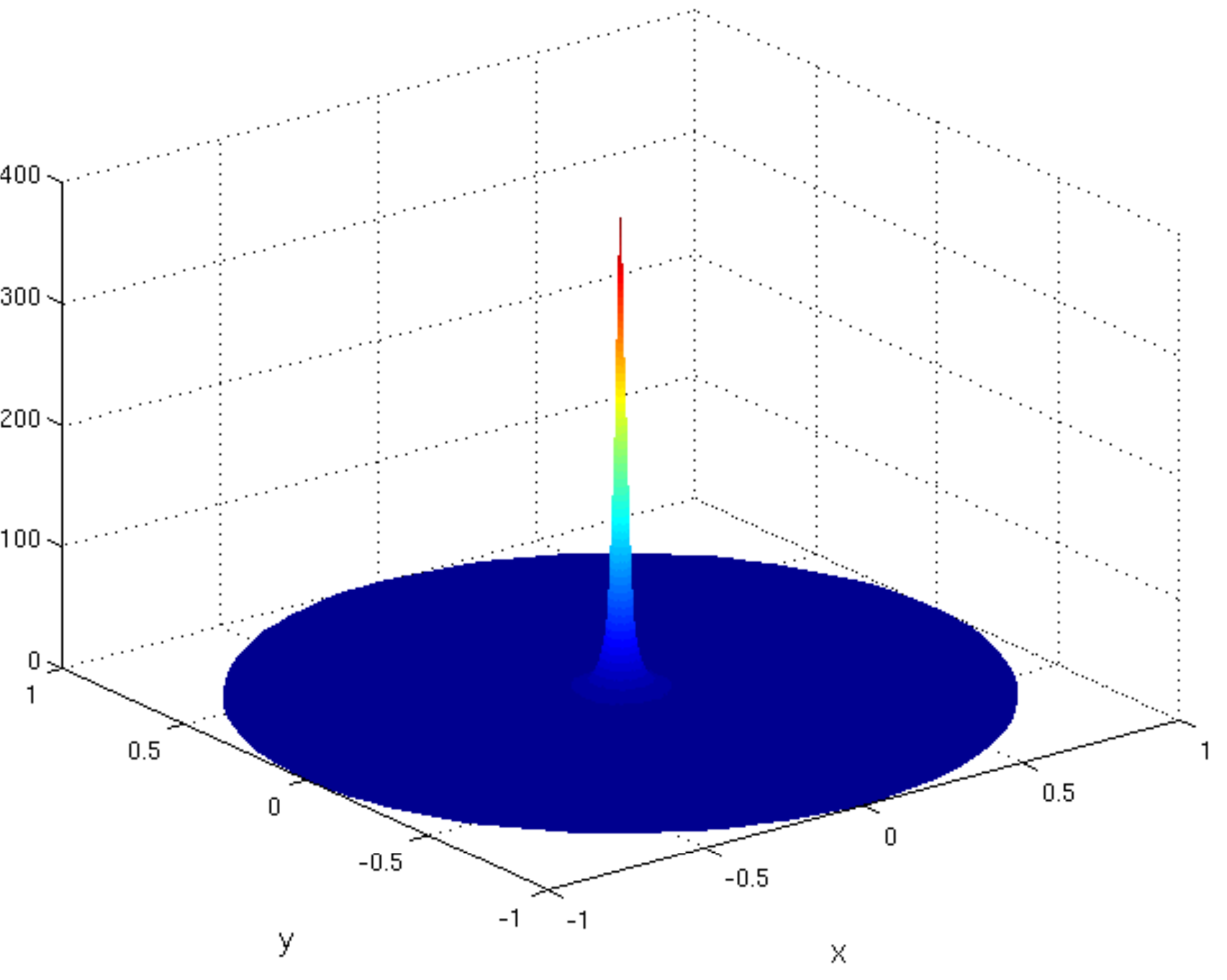}}
\subfloat[Relative entropy]{\includegraphics[width=0.32\textwidth]{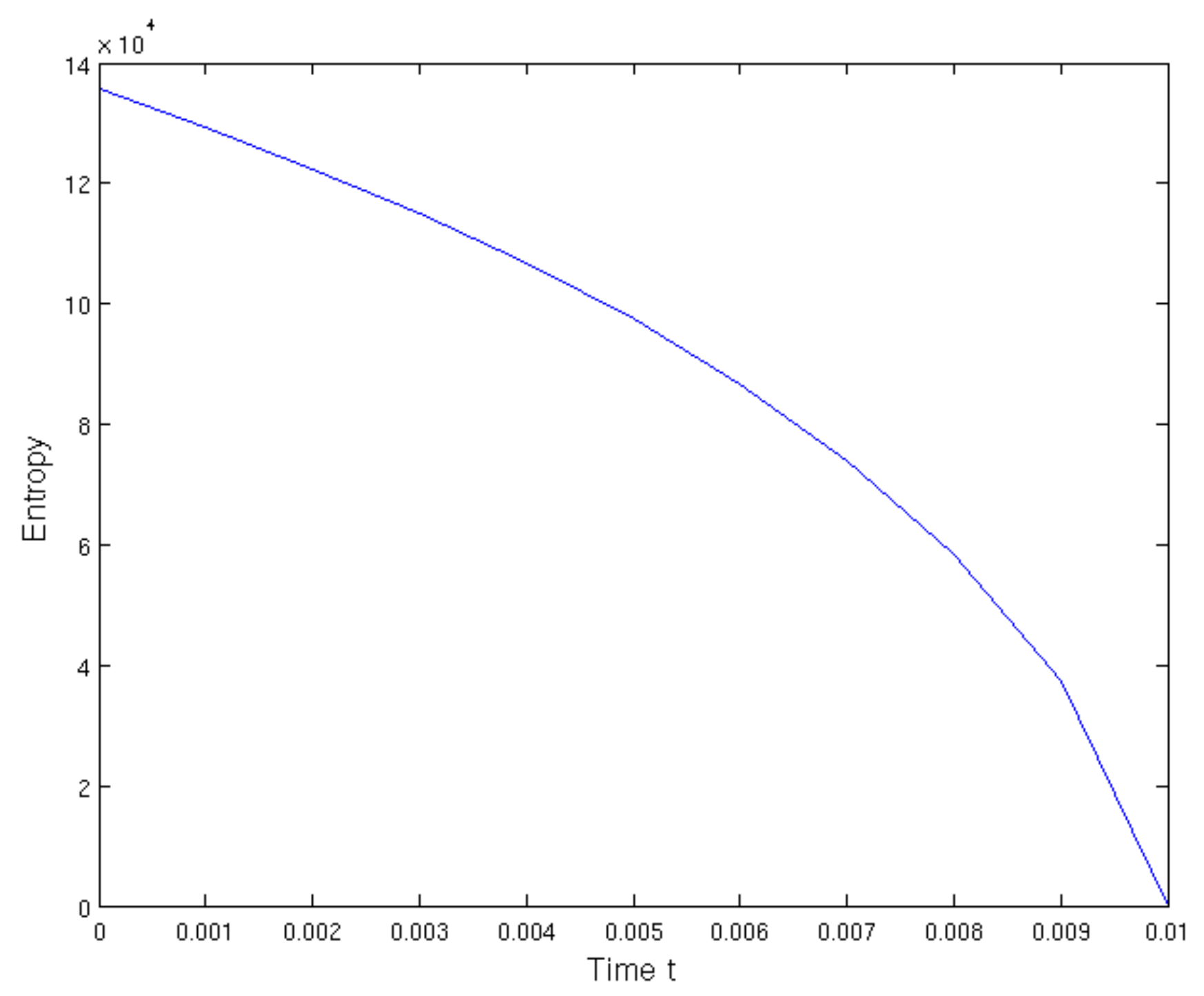}}
\caption{Simulation results for KS model with initial $M=1$ and $\chi = 1.1 \times 8\pi$.} \label{f:ex5}
\end{center}
\end{figure}

In our final example we consider the modified KS model \eqref{e:mksm} with an initial Gaussian of mass one and $\chi = 1.1 \times 8 \pi$. The time steps are set to $\Delta t=5\times 10^{-3}$ and we observe the fast formation of a Delta Dirac at the center of the domain, see Figure \ref{f:ex5}. Figure \ref{f:ex5}(c) indicates that the free energy decay is changing concavity as the free energy tries to decay faster possibly tending to $-\infty$ at the blow-up time.

\section{Conclusion}\label{s:conclusion}

\noindent In this paper we propose a numerical algorithm for nonlinear PDEs, which can be written as gradient flows with respect to the Wasserstein distance. The method
is based on its gradient flow formulation with respect to the Euclidean distance, corresponding to the Lagrangian representation of the original nonlinear PDE. 
The construction of the solver guarantees the preservation of structural features, such as entropy decay or positivity of solutions. We presented extensive numerical simulations 1D, which
illustrated the flexibility of our approach and confirmed theoretical results concerning entropy decay and blow up behavior. Even though these numerical simulations confirmed the predicted 
convergence behavior towards the steady state, the numerical analysis is still an open problem (even in 1D). Hence we  plan to investigate the numerical analysis of the proposed scheme in spatial dimension one and as a next step study the method for radially symmetric solutions.

The 2D algorithm involves several pre- and post-processing steps, which present additional challenges with respect to numerical accuracy and computational complexity. The proposed
finite element discretization allows to consider more general domains and better resolve features of radially symmetric solutions, which often arise in aggregation equations. However the
computational complexity of the preprocessing step as well as the solver itself remains a key limitation of the solver, which we plan to address in a future work.

\section*{Acknowledgments}
\noindent JAC was partially supported by the Royal Society via a Wolfson Research Merit Award. HR and MTW acknowledge financial support from the Austrian Academy of Sciences \"OAW via the New Frontiers Group NSP-001. The authors would like to thank the King Abdullah University of Science and Technology for its hospitality and partial support while preparing the manuscript.

\section*{References}
\bibliographystyle{plain} 
\bibliography{carrilloranetbauerwolfram}

\end{document}